\documentclass[12pt]{amsart}
\usepackage{amsmath,amssymb,amsfonts,amscd,graphicx}
\usepackage{xypic}
\newtheorem{theorem}[subsection]{Theorem}
\newtheorem{lemma}[subsection]{Lemma}

\newtheorem{proposition}[subsection]{Proposition}

\theoremstyle{definition}
\newtheorem{definition}[subsection]{Definition}

\theoremstyle{remark}

\numberwithin{equation}{section}

\theoremstyle{definition}

\newcommand{\co}{\colon\thinspace}

\newcommand{\nl}{\hfil\break}

\begin{document}

\title{Topological arbiters}
\author[Michael Freedman and Vyacheslav Krushkal]{Michael Freedman and Vyacheslav Krushkal$^{^*}$}
\address{Microsoft Station Q, University of California, Santa Barbara, CA 93106-6105}
\email{michaelf\char 64 microsoft.com}
\address{Department of Mathematics, University of Virginia, Charlottesville, VA 22904}
\email{krushkal\char 64 virginia.edu}

\thanks{$^{^*}$Partially supported by the NSF}

\begin{abstract}
This paper initiates the study of topological arbiters, a concept rooted in
Poincar\'{e}-Lefschetz duality. Given an
$n$-dimensional manifold $W$, a topological arbiter associates a value $0$ or $1$ to
codimension zero submanifolds of $W$, subject to natural topological and duality axioms.
For example, there is a unique arbiter on ${\mathbb R}P^2$, which reports the location of the
essential $1$-cycle.
In contrast, we show that there exists an uncountable collection of topological arbiters in dimension $4$.
Families of arbiters, not induced by homology, are also shown to exist in higher dimensions.
The technical ingredients underlying the four dimensional results are secondary obstructions to generalized
link-slicing problems. For classical links in $S^3$ the
construction relies on the existence of nilpotent embedding obstructions in dimension $4$, reflected
in particular by the Milnor group. In higher dimensions novel arbiters
are produced using nontrivial squares in stable homotopy theory.

The concept of ``topological arbiter'' derives from percolation and from $4$-dimen- sional surgery. It is not the purpose of this paper to advance either of these subjects, but rather to study the concept for its own sake. However in appendices we give both an application to percolation, and the current understanding of the relationship between arbiters and surgery.
An appendix also introduces a more general notion of a multi-arbiter.
Properties and applications are discussed, including a construction of non-homological multi-arbiters.
\end{abstract}

\maketitle

\section{Introduction.} Given a closed smooth $n-$manifold $W$,
consider the collection ${\mathcal M}_W$ of connected smooth codimension zero
submanifolds of $W$.

\begin{definition} \label{def: arbiter}
A {\em topological arbiter} on $W$ is a function ${\mathcal A}\co
{\mathcal{M}}_W\longrightarrow \{0,1\}$ satisfying axioms (1) -- (3):

(1) \parbox[t]{14.25cm}{{\em ``${\mathcal A}$ is topological''}: If $M, M'\in {\mathcal M}_W$ and $M$ is ambiently isotopic
to $M'$ in $W$ then ${\mathcal A}(M)={\mathcal A}(M')$.}

(2) \parbox[t]{14.25cm}{{\em ``Greedy axiom''}: If $M\subset M'$
and ${\mathcal A}(M)=1$ then ${\mathcal A}(M')=1$.}

(3) \parbox[t]{14.25cm}{{\em ``Duality''}: Suppose $A,B\in{\mathcal M}_W$ are such that $W=A\cup B$, with $A\cap B
=\partial A\cap \partial B$.
Then ${\mathcal
A}(A)+{\mathcal A}(B)=1$.}
\end{definition}

The notion of a topological arbiter was introduced by the first author in \cite{F0} (the term {\em exclusive property} was used
in that reference) in connection with percolation theory. The purpose of this paper is to investigate further the structure of these
invariants in $4$ and higher dimensions.

A prototypical example of an arbiter is the homological arbiter ${\mathcal A}_h$ on ${\mathbb R}P^2$. Specifically,
given $M\subset {\mathbb R}P^2$ define ${\mathcal A}_h(M)=1$ if $M$ carries the non-trivial
first homology class of ${\mathbb R}P^2$, i.e. if $$H_1(M; {\mathbb Z}/2)\longrightarrow
H_1({\mathbb R}P^2; {\mathbb Z}/2)\cong {\mathbb Z}/2$$ is onto, and set ${\mathcal A}_h(M)=0$ otherwise. The first two axioms
are immediate, and (3) easily follows from Poincar\'{e} duality, see lemma \ref{projective space} below.
In fact, ${\mathcal A}_h$ is the {\em unique} topological arbiter on the projective plane. Indeed, suppose there is a different
arbiter ${\mathcal A}$ on ${\mathbb R}P^2$, so there is $A\subset {\mathbb R}P^2$ such that
${\mathcal A}(A)=1$ but  $H_1(A; {\mathbb Z}/2)\longrightarrow
H_1({\mathbb R}P^2; {\mathbb Z}/2)$ is the trivial map. Consider the complement
$B={\mathbb R}P^2\smallsetminus A$, then by Poincar\'{e}-Lefschetz duality $B$ carries the first homology of ${\mathbb R}P^2$.
Since $A$ lies in the complement of a non-trivial cycle, it is contained in a $2$-cell $D^2\subset {\mathbb R}P^2$. It follows from the axioms (2), (3) that ${\mathcal A}(D^2)=0$, so by (2) ${\mathcal A}(A)=0$, a contradiction proving that ${\mathcal A}_h$ is the unique
arbiter on the projective plane.

It is easy to see that for any $n$, the $n$-sphere $S^n$ does not admit a topological arbiter:
consider a decomposition $S^n=A\cup B$ where both $A,B$ are diffeomorphic to $D^n$. According to axiom (1), ${\mathcal A}(A),
{\mathcal A}(B)$ are either both $0$ or both $1$, contradicting the duality axiom (3).

More generally, note that any manifold $W$ admitting an open book
decomposition does not support a topological arbiter. (However $W$ may admit a {\em multi-}arbiter, see the Appendix.)
Indeed, in this case the complement
of the binding fibers over $S^1$, where the binding is a specified codimension $2$ submanifold of $W$.
Decompose $S^1$ into a union of two intervals and consider the corresponding decomposition
$W=A\cup B$. $A$ and $B$ are diffeomorphic (in fact, isotopic), so again axioms (1), (3)
contradict the existence of an arbiter. Note that odd-dimensional manifolds admit open book decompositions \cite{Lawson, Winkelnkemper}.
Simply-connected even-dimensional
manifolds (of dimension $\geq 6$) are open books, provided that the signature is zero \cite{Winkelnkemper}.
In the general even dimensional case there is a further obstruction given by nonsingular bilinear forms \cite{Quinn}.

Generalizing the discussion of the projective plane above, consider the ``homological'' arbiter on even-dimensional
real and complex projective spaces:

\begin{lemma} \label{projective space} \nl \sl
(1) Let $M\subset{\mathbb R}P^{2n}$.
Define ${\mathcal A}_h(M)=0$ if $H_n(M;{\mathbb Z}/2)\longrightarrow
H_n({\mathbb R}P^{2n}; {\mathbb Z}/2)\cong {\mathbb Z}/2$ is the trivial map and
${\mathcal A}_h(M)=1$ otherwise.

(2) Fix a field $F$ and let $M\subset{\mathbb C}P^{n}$, where $n$ is even. Define ${\mathcal
A}_F(M)=0$ if $H_n(M;F)\longrightarrow H_n({\mathbb C}P^{n}; F)\cong F$ is
the trivial map and ${\mathcal A}_F(M)=1$ otherwise.

Then ${\mathcal A}_h$ is a topological arbiter on
${\mathbb R}P^{2n}$, and for any $F$, ${\mathcal A}_F$ is a
topological arbiter on ${\mathbb C}P^{n}$.
\end{lemma}

{\em Proof}. The first two axioms are immediate. To prove axiom (3), let $W={\mathbb R}P^{2n}$ or ${\mathbb
C}P^{n}$, let the coefficients $F={\mathbb Z}/2$ in the first case
and an arbitrary field in the second case. Consider the long exact
sequences
\begin{equation} \label{comm diagram}
\xymatrix{
    H_n(A) \ar[r]
    \ar[d] &
    H_n(W)
    \ar[d] \ar[r] &
    H_n(W,A) \ar[d]
    \\
    H^n(W,B) \ar[r]
    &
    H^n(W)
    \ar[r] &
    H^n(B)
    }
\end{equation}
where the vertical maps are isomorphisms given by
Poincar\'{e}-Lefschetz duality (cf \cite{Hatcher}). Since $H_n(W;F)\cong F$ it is clear
that precisely one of the two maps $H_n(A;F)\longrightarrow H_n(W;F)$,
$H_n(B;F)\longrightarrow H_n(W;F)$ is non-trivial. \qed

There are further examples of manifolds and arbiters induced by homology.
For example, generalizing lemma \ref{projective space} one may define homological arbiters on any
even-dimensional manifold $W^{2n}$ whose middle-dimensional homology group $H_n(W;F)$ has odd rank.
The arbiter ${\mathcal A}_F$ associates $1$ to submanifolds representing more than half the rank of
$H_n(W;F)$.

The results in this paper concern arbiters on even-dimensional balls $D^{2n}$, the local version
of topological arbiters discussed above. {\em Local arbiters} (see definition \ref{def: local arbiter}
in section \ref{sec: Local arbiters} below) associate a value $0$ or $1$ to codimension zero submanifolds of $D^{2n}$
with a prescribed boundary condition: $M\cap \partial D^{2n}$ is a regular neighborhood of the standard
$(n-1)$-sphere $S^{n-1}$ in $S^{2n-1}=\partial D^{2n}$, $\partial M\cap \partial D^{2n}
=\partial D^n\times D^n$, writing $D^{2n}=D^n\times D^n$. The duality axiom satisfied by local arbiters is modeled on
the Alexander duality satisfied by submanifolds $D^{2n}=A\cup B$ bounding the Hopf link $S^{n-1}\cup S^{n-1}
\subset S^{2n-1}$. The reader may wish to glance ahead to absorb definition \ref{def: local arbiter}.

All arbiters discussed so far are induced by homology (with various coefficients).
It follows from the universal coefficient theorem that
the arbiter ${\mathcal A}_F$ on ${\mathbb C}P^n$ (and the analogous local arbiter on $D^{2n}$)
depends only on the characteristic of the
field $F$, so $\{ {\mathcal A}_F \}$ is a countable collection.
Focusing on $D^4$, we construct
a collection of arbiters different from the homological ones:

\begin{theorem} \label{uncountable} \sl There exists an uncountable collection of local topological arbiters in dimension $4$.
\end{theorem}

The construction (discussed in more detail below and in section \ref{sec: Local arbiters})
relies on the existence of non-abelian embedding obstructions in $D^4$, reflected
in particular by the Milnor group. These obstructions, detecting in  a certain sense the failure of the Whitney trick in $4$ dimensions,
have been extensively studied in the setting of link-homotopy theory \cite{Milnor}.

Since the Whitney trick is valid in dimensions $> 4$, our construction of (uncountably
many) arbiters on $D^4$ does not have an immediate analogue in higher dimensions.
In section \ref{sec: Higher dimensions} we give a rather different construction of families of local arbiters, not induced by homology, in
dimensions $6$ and higher. These higher-dimensional arbiters
are homotopy-theoretic in nature, the following theorem summarizes our construction:

\begin{theorem} \label{higher dim thm}\sl
To each non-trivial square in the stable homotopy ring of spheres
there is associated a local topological arbiter not induced by homology on
$D^{2n}$ for sufficiently large $n$. In particular, there exist local arbiters not induced by homology on $D^{2n}$ for each
$n>2$.
\end{theorem}

The main ingredient used in the construction of these arbiters in higher dimensions is the following homotopy-theoretic obstruction
to a generalized link-slicing problem: for each $n>2$ we show that there exist submanifolds $(M,S^{n-1})\subset (D^{2n},S^{2n-1})$,
where $S^{n-1}$ is a distinguished ``attaching'' sphere in the boundary of $M$,  such that

(1) \parbox[t]{14.25cm}{for any coefficient ring $R$ the attaching $(n-1)$-sphere
of $M$ is non-trivial in $H_{n-1}(M; R)$, and}

(2) \parbox[t]{14.25cm}{The components of the Hopf link $S^{n-1}\cup S^{n-1}\subset S^{2n-1}$ do not bound
disjoint embeddings of two copies of $M$ into $D^{2n}$.}

Observe that if the attaching sphere $S^{n-1}$ were trivial in $H_{n-1}(M;F)$ for some field coefficients $F$, then Alexander
duality would immediately imply (2). The proof of (1) and (2), given in section \ref{sec: Higher dimensions},
can be carried out using a ``secondary'' obstruction to disjointness where
$n+1$-cells are attached to $S^{n-1}$ via the generator of ${\pi}_1^S$, $n>3$. A further refinement is necessary for $D^6$, see section \ref{D6}.

To prove theorem \ref{higher dim thm}, one defines a ``partial'' arbiter, setting it equal to $1$ on all submanifolds of $D^{2n}$ containing
$M$ and equal to zero on all submanifolds contained in $D^{2n}\smallsetminus M$. The condition (2) above implies that this partial arbiter
is consistently defined, and (1) shows that it is different from any homological arbiter. The proof of the theorem is completed
by extending this partial arbiter to all submanifolds of $D^{2n}$, the details are given in section \ref{sec: Higher dimensions}.

A brief outline of the rest of the paper is as follows.
The proof of theorem \ref{uncountable} begins in section \ref{sec: Local arbiters}
with the introduction of a dyadic collection of ``model'' submanifolds of $D^4$ (i.e.
a special collection of submanifolds, parametrized by vertices of a trivalent tree.)
A generalization of the
Milnor group is used to find obstructions to pairwise embeddings of these model submanifolds. The collection ${\mathcal M}_{D^4}$ of all submanifolds
of $D^4$ has a partial ordering induced by inclusion.
The non-embedding results mentioned above allow one to
define an uncountable collection of arbiters on those elements of ${\mathcal M}_{D^4}$ which are comparable
to model submanifolds. One then extends these arbiters to ${\mathcal M}_{D^4}$, concluding the proof of theorem
\ref{uncountable}, see section \ref{sec:proof}.

Section \ref{sec: embeddings} points out that the embeddings of the submanifolds $M\hookrightarrow W$ in definition \ref{def: arbiter}
of arbiters play a very important role. The result of \cite{K0} shows that if one considers topological arbiters as
invariants of manifolds, disregarding their embedding information, then in the context of $D^4$ whenever the attaching circle
of $M\subset D^4$ is trivial in $H_1(M;{\mathbb Z})$ the value of any arbiter on $M$ has to be $1$.
The proof of theorem \ref{uncountable} is valid only when the embedding data is taken into account. Indeed
the arbiters constructed in its proof assume value $1$ on manifolds where the attaching circle does not vanish in homology with any
coefficients.

Section \ref{sec: CP2} adapts
the local construction in the proof of theorem \ref{uncountable}
to the setting of ${\mathbb C}P^2$. It seems likely that this construction should yield an
uncountable collection of arbiters on ${\mathbb C}P^2$, however proving this would require a substantial
extension of the techniques of Milnor's theory of link homotopy.

In appendix \ref{sec: Surgery} we note that the existence of a local arbiter on $D^4$ satisfying an
additional ``Bing doubling'' axiom would give an obstruction to topological $4$-dimensional surgery for
the non-abelian free groups, using the A-B slice reformulation \cite{F, F1} of the surgery conjecture.

Appendix \ref{sec: multiarbiters} introduces two generalizations of topological arbiters, {\em multi-arbiters}.
We discuss their properties and applications, in particular a construction of non-homological multi-arbiters,
and an application to percolation theory.

\section{Local arbiters.} \label{sec: Local arbiters}

The purpose of this section is to prove theorem \ref{uncountable}.
We start by
defining the notion of a local topological arbiter on $D^4$ (more generally, the definition
makes sense for $D^{2n}$ for any even dimension $2n$.)
Consider
$\mathcal{M}=\mathcal{M}_{D^4}=\{(M,{\gamma})|\, M$ is a codimension zero, smooth,
compact submanifold of $D^4$, ${\gamma}\subset
\partial M$, and $M\cap\partial D^4$ is a tubular neighborhood of an
unknotted circle ${\gamma}\subset S^3 \}$. (For $n>2$, $M\cap \partial D^{2n}$
is a tubular neighborhood of a standard $S^{n-1}$, isotopic in $\partial D^{2n}
=\partial (D^n\times D^n)$ to $\partial D^n\times D^n$.)
Consider a local version of definition \ref{def: arbiter}:

\begin{definition} \label{def: local arbiter}
A {\em local topological arbiter} is a function ${\mathcal A}\co
{\mathcal{M}}\longrightarrow \{0,1\}$ satisfying axioms (1) -- (3):

(1) \parbox[t]{14.25cm}{If $(M,{\gamma}), (M',{\gamma}')\in {\mathcal M}$ and $(M,{\gamma})$ is ambiently isotopic
(not necessarily relative to the identity on the boundary) to $(M',{\gamma}')$ in $(D^4,S^3)$ then ${\mathcal A}(M,{\gamma})={\mathcal A}(M',{\gamma}')$.}

(2) \parbox[t]{14.25cm}{If $(M,{\gamma})\subset (M',{\gamma}')$ and
${\mathcal A}(M,{\gamma})=1$ then ${\mathcal A}(M',{\gamma}')=1$.}

(3) \parbox[t]{14.25cm}{Suppose $(A,{\alpha}),(B,{\beta})\in{\mathcal M}_{D^4}$ are such that $D^4=A\cup B$, with $A\cap B
=\partial A\cap \partial B$, and so that the attaching curves ${\alpha}, {\beta}$ of $A,B$
form the Hopf link in $\partial D^4$. (In other words, $D^4=A\cup B$ is a decomposition of
$D^4$ extending the standard genus $1$ Heegaard decomposition of the
$3-$sphere.) Then ${\mathcal
A}(A,{\alpha})+{\mathcal A}(B,{\beta})=1$. (For $n>2$ the ``Hopf link'' above is replaced with
$\partial D^n\times 0\sqcup 0\times \partial D^n\subset S^{2n-1}.$)}
\end{definition}

We summarize local versions of the properties discussed in the introduction.
In dimension $4$ (more generally, in any even dimension), given a field $F$ there is a ``homological''
local arbiter ${\mathcal A}_F$. This arbiter (by the universal coefficient theorem depending
only on the characteristic of $F$) is defined by ${\mathcal A}_F(M,{\gamma})=1$ if and only if ${\gamma}=0\in
H_1(M;F)$. The axioms (1), (2) are satisfied automatically, (3) is a
consequence of Alexander duality.

In dimension $2$ there is a unique local arbiter ${\mathcal A}$ (equal to ${\mathcal A}_F$ for any
field $F$.)
The difference between arbiters in dimensions $2$ and higher is explained by the divergence
between homology and homotopy: for $(A,{\alpha})\subset (D^2,S^1)$, ${\alpha}=0\in H_0(A)$ iff ${\alpha}=0\in {\pi}_0(A)$.
If this is the case, the complement $B$ of $A$ in $D^2$ may be pushed off an arc connecting $\alpha$ in $A$ into a collar
on its attaching curve $\beta$, therefore the value of any arbiter on $B$ is zero.
In higher dimensions, being trivial in homology is a weaker condition, and we show that there exist arbiters
whose value on $B$ may be $1$ even though $\beta$ is homologically non-trivial with any coefficients.

\subsection{Model decompositions.} \label{sec: Model decompositions}

The proof of theorem \ref{uncountable} starts with the construction of a special collection of submanifolds of $D^4$, which are shown
below to satisfy pairwise non-embedding properties.  These codimension zero
submanifolds form complementary pairs $A_I, B_I$ parametrized by a dyadic multiindex $I$ discussed in more
detail below. For each index $I$,
the pair is a decomposition of the $4$-ball in the sense of definition \ref{def: local arbiter} (3) above: $D^4=A_I\cup B_I$, and
the distinguished curves ${\alpha}_I, {\beta}_I$ form the Hopf link in $\partial D^4$.
To start the construction, consider the genus one surface
$S$ with a single boundary component, and set $B_1=S\times D^2$.
Denote ${\beta}_1=\partial S\times 0$.
To specify an embedding $(B_1,{\beta}_1)\subset (D^4, S^3)$, embed the surface $S$ into $S^3$, push
it into the $4$-ball and take a regular neighborhood.

\begin{figure}[ht]
\includegraphics[width=2.6cm]{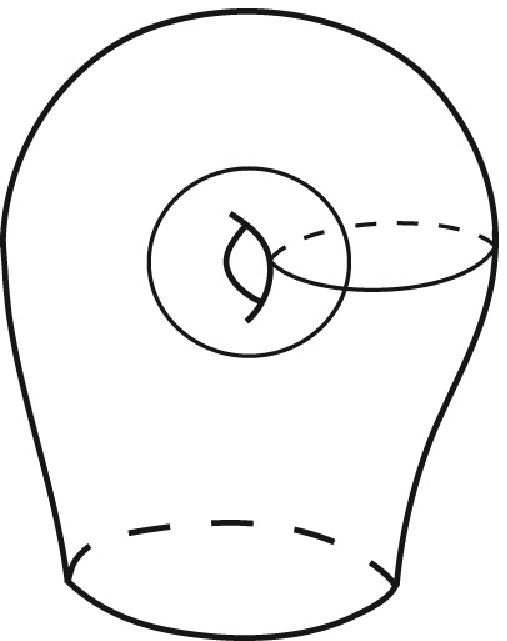} \hspace{3cm} \includegraphics[width=3.5cm]{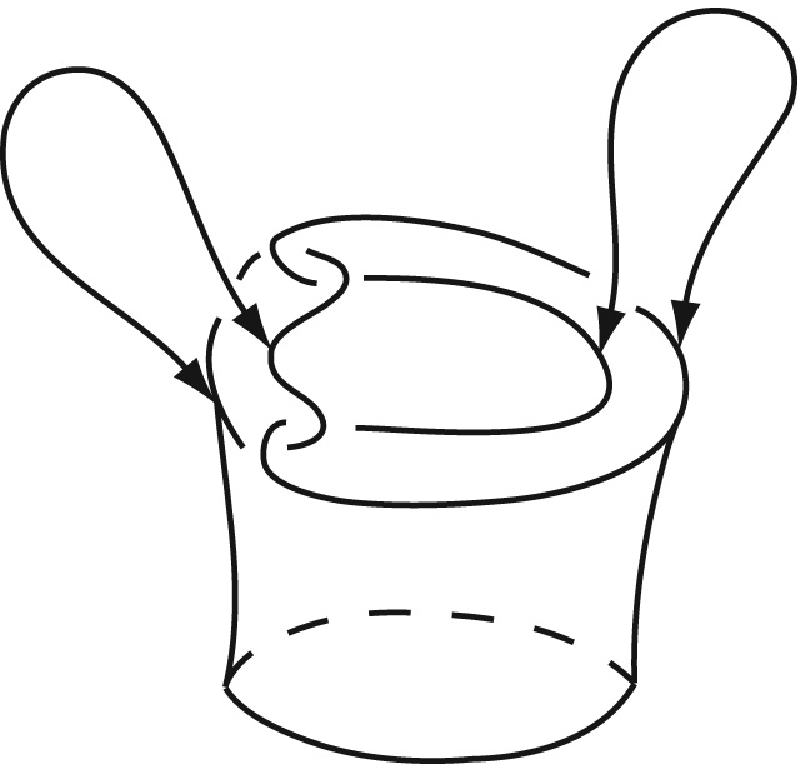}
{\small
    \put(-290,37){$B_1$}
   \put(-258,48){${\gamma}_1$}
    \put(-212,42){${\gamma}_2$}
    \put(-270,-6){${\beta}_1$}
    \put(-80,-6){${\alpha}_1$}
    \put(-116,73){$H_1$}
    \put(3,80){$H_2$}
    \put(-8,23){$A_1$}}
\caption{}
\label{fig:surface Bing double}
\end{figure}

\begin{figure}[ht]
\vspace{.3cm}
\includegraphics[width=4cm]{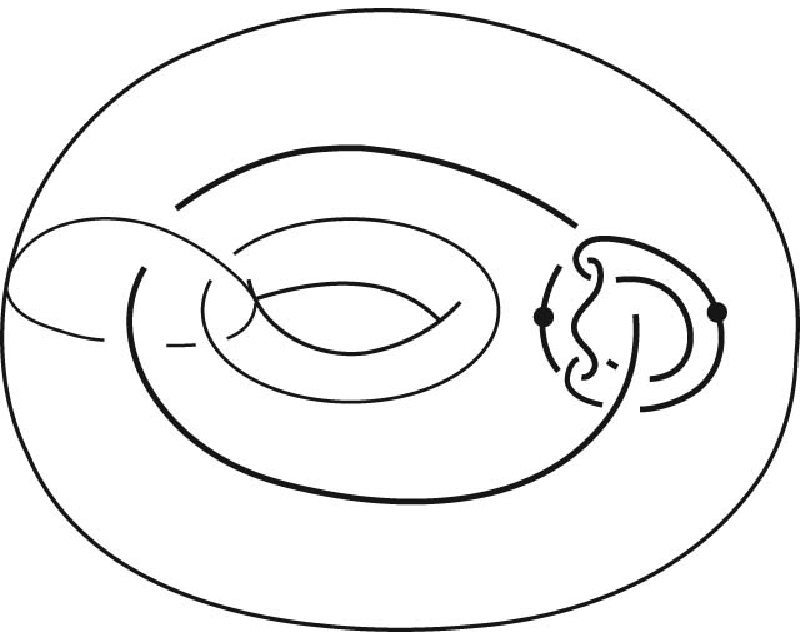} \hspace{3cm} \includegraphics[width=4cm]{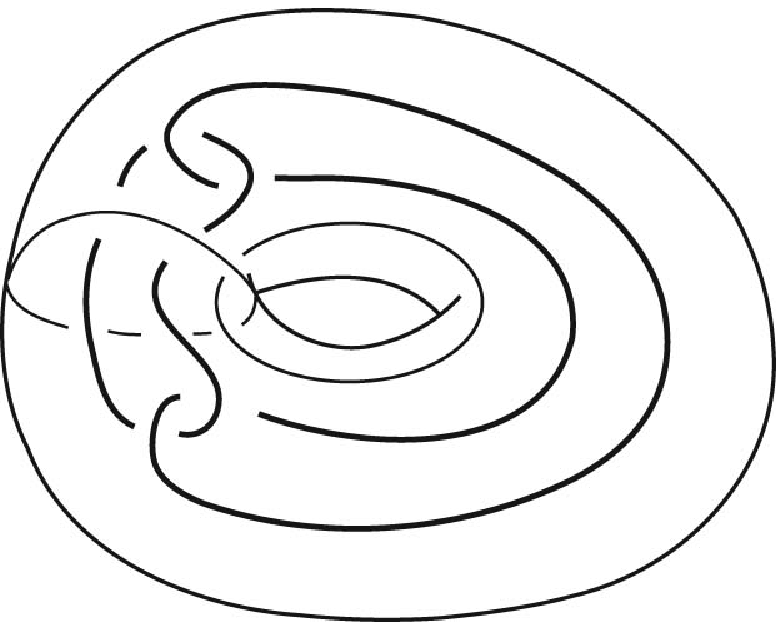}
{\small
    \put(-336,5){$B_1$}
    \put(-128,5){$A_1$}}
{\scriptsize
    \put(-304,20){$0$}
    \put(-259,28){${\beta}_1$}
    \put(-107,30){$0$}
    \put(-17,22){$0$}
    \put(-46,35){${\alpha}_1$}}
\caption{}
\label{fig:Kirby Bing double}
\end{figure}

\begin{proposition} \label{surface complement}
\sl The complement $A_1=D^4\smallsetminus B_1$ is obtained from the collar
$S^1\times D^2\times I$ by attaching a pair
of zero-framed $2$-handles to the Bing double of the core of the
solid torus $S^1\times D^2\times 1$, figures \ref{fig:surface Bing double}, \ref{fig:Kirby Bing double}.
\end{proposition}

 The distinguished curve ${\alpha}_1\subset \partial A_1$ is formed by $S^1\times 0 \times 0$.
Figure \ref{fig:surface Bing double} is a schematic ``spine'' picture of
the manifolds $A_1, B_1$, a precise description in
terms of Kirby diagrams (drawn in a solid torus) is given in figure \ref{fig:Kirby Bing double}.
Note that a distinguished pair of curves ${\gamma}_1, {\gamma}_2$,
forming a symplectic basis in the surface $S$, is determined as the
meridians (linking circles) to the cores of the $2-$handles $H_1,
H_2$ in $D^4$. In other words, ${\gamma}_1$,
${\gamma}_2$ are fibers of the circle normal bundles over the cores
of $H_1, H_2$ in $D^4$.
The proof of proposition \ref{surface complement} is an exercise in Kirby calculus, see for example
\cite{FL}.

The manifolds $A_1, B_1$ are the basic building blocks in the inductive construction of model decompositions.
Consider the $2$-handle $H_1$ (one of the $2$-handles in the definition of $A_1$, figure \ref{fig:surface Bing double})
in place of the original
$4$-ball. The pair of curves (${\gamma}_1$, the attaching circle
${\delta}_1$ of $H_1$) forms the Hopf link in the boundary of $H_1$.
In $H_1$ consider the standard genus one surface bounded by
${\delta}_1$. As discussed above, its complement is given by two
zero-framed $2$-handles attached to the Bing double of ${\gamma}_1$.
Assembling this data, consider the new decomposition $D^4=A_0\cup
B_0$. The manifold $A_0$ and its Kirby diagram are shown in figure \ref{fig:A0} below, and both
$A_0, B_0$ are included in figure \ref{fig:model decompositions}.

We are now in a position to define the more general collection of special
elements $A_I, B_I$ of ${\mathcal M}$ (``model decompositions'',
$D^4=A_I\cup B_I$), figure \ref{fig:model decompositions}. These decompositions correspond to
vertices of a trivalent tree, and each $A_I, B_I$ is labeled by a
dyadic multiindex $I$. Each arrow in figure \ref{fig:model decompositions} represents a
an inclusion. These model decompositions are defined
inductively; the transition from $A_I$ to $A_{I,0}$ replaces one of
the two $2$-handles of $A_I$, introduced during the previous step,
by $A_0$. $A_{I,1}$ is obtained by replacing one of these two
``new'' $2$-handles of $A_I$ by $A_1$ (in other words, by Bing
doubling this $2$-handle.) Then $B_I$ is defined as the complement
$D^4\smallsetminus A_I$.  Note that this is a subset of the collection of
decompositions of $D^4$ introduced in the context of the A-B slice problem in \cite{FL}.
Figure \ref{fig:model decompositions} shows only the ``spines'' of the submanifolds
$A_I, B_I$; we refer the reader to \cite{FL} for a precise description of
these $4$-manifolds in terms of
Kirby diagrams.

\begin{figure}[ht]
\includegraphics[height=8cm]{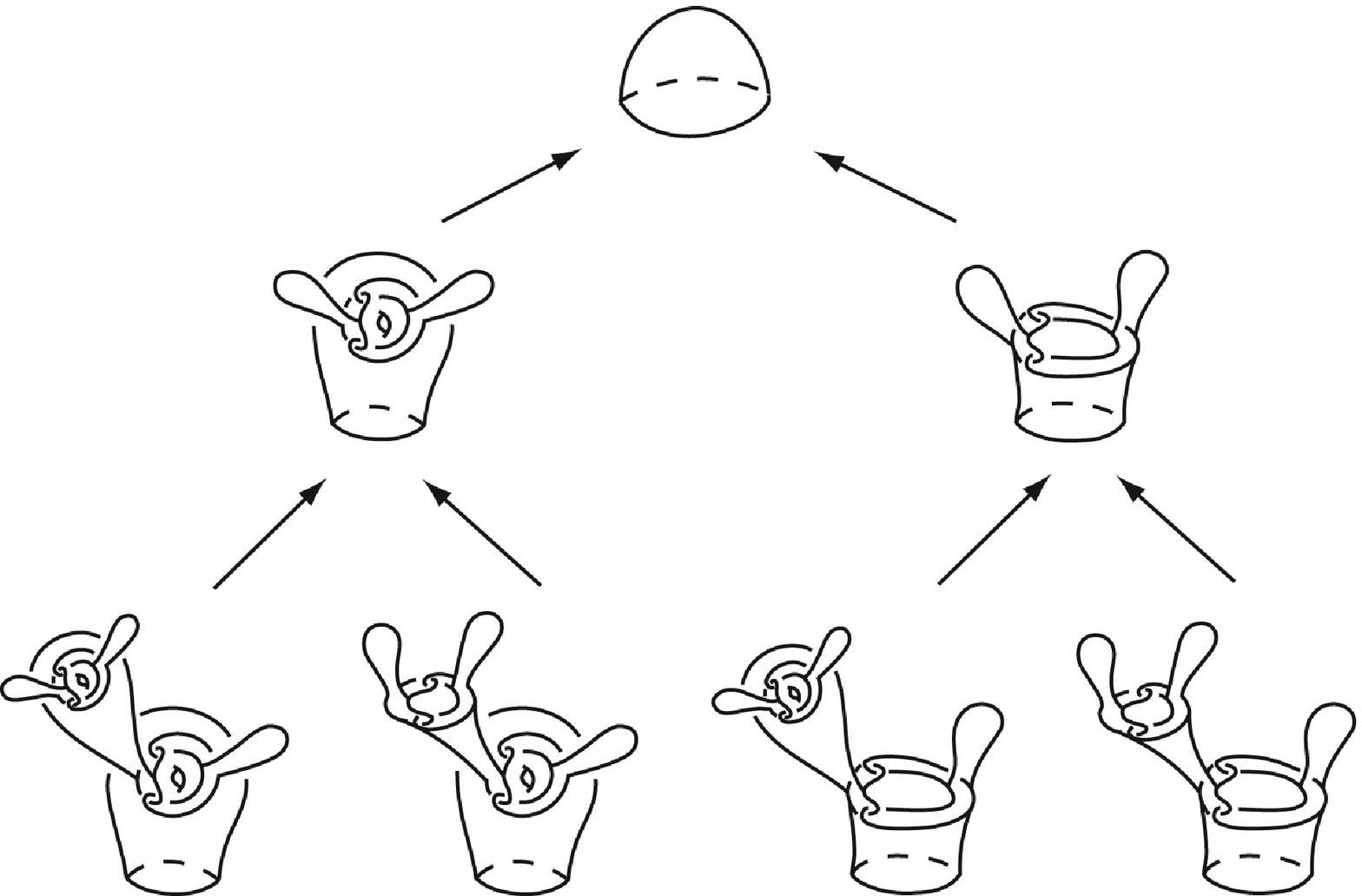}
{\footnotesize
    \put(-200,209){$A$}
    \put(-280,127){$A_0$}
    \put(-103,127){$A_1$}
    \put(-340,12){$A_{00}$}
    \put(-250,12){$A_{01}$}
    \put(-150,12){$A_{10}$}
    \put(-58,12){$A_{11}$}
    \put(-180,185){${\alpha}$}
    \put(-238,113){${\alpha}_0$}
    \put(-61,113){${\alpha}_1$}}
\includegraphics[height=8cm]{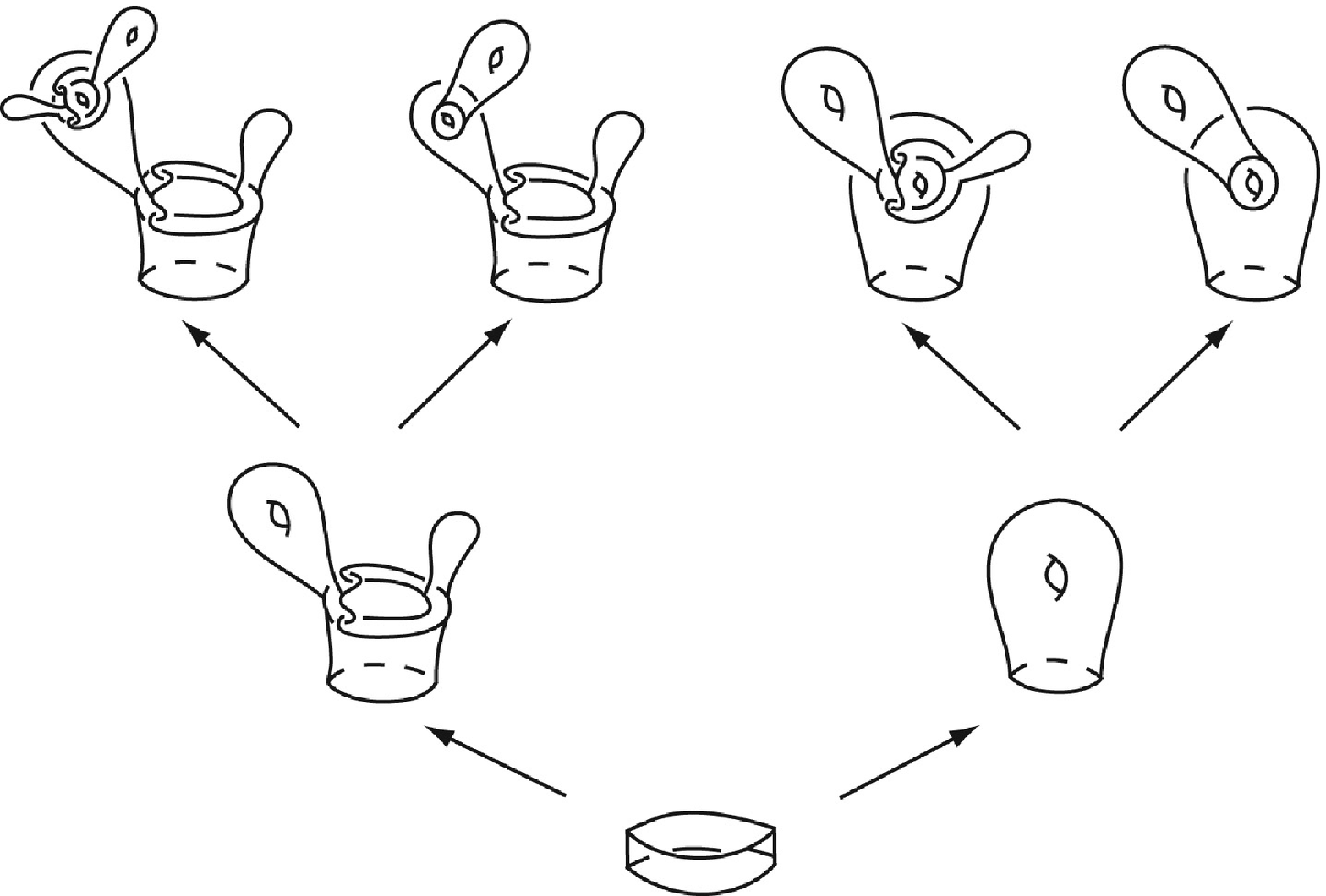}
{\footnotesize
    \put(0,265){}
    \put(-190,7){$B$}
    \put(-274,63){$B_0$}
    \put(-99,65){$B_1$}
    \put(-323,165){$B_{00}$}
    \put(-233,165){$B_{01}$}
    \put(-140,165){$B_{10}$}
    \put(-53,165){$B_{11}$}
    \put(-146,-6){${\beta}$}
    \put(-223,45){${\beta}_0$}
    \put(-53,43){${\beta}_1$}}
\caption{The collection of model decompositions $A_I, B_I$. The
arrows represent inclusions.}
\label{fig:model decompositions}
\end{figure}

Consider the {\em lexicographic order} on the set of dyadic
multiindices: given $I=(i_1,\ldots,i_m)$, $J=(j_1,\ldots,j_n)$, set
$I<J$ if $i_1=j_1,\ldots, i_k=j_k$ and $i_{k+1}=0, j_{k+1}=1$ for
some $k$. Also set $I<J$ if $I$ is an initial segment of $J$:
$J=(I,j_{k+1},\ldots$)

Consider the trivalent tree $T$ encoding the submanifolds $A_I$. The root
of $T$ corresponds to $A=2$-handle. Let $r$ be a ``geodesic'' ray in $T$
originating at the root; $r$ is viewed as an infinite dyadic multiindex.
For a technical reason which will be apparent in the proof below, it is convenient to consider all rays $r$ passing
through $A_1$ (i.e. whose first coefficient is $1$),
and the lexicographic order is extended
to the set of such rays.

\begin{definition}
Let ${\mathcal M}_0\subset {\mathcal M}$ denote the collection of submanifolds
$\{ A_I, B_I\}$ defined above.
For each $r$, a ray as above, define a {\em partial arbiter}
${\mathcal A}_r\co{\mathcal M}_0\longrightarrow \{0,1\}$ by setting
$${\mathcal A}_r(A_I)=1 \,\, {\rm if} \, \, I<r, \,\, {\rm and}\, \,
{\mathcal A}_r(A_I)=0\, \, {\rm if}\,\,  I>r.$$ Then the
duality axiom $(3)$ in definition \ref{def: local arbiter} forces one to define $${\mathcal A}_r(B_I)=0
\,\, {\rm if}\,\, I<r, \,\, {\rm and} \,\, {\mathcal A}_r(B_I)=1\,\,
{\rm if} \,\, I>r.$$
\end{definition}

\medskip

\begin{lemma} \label{greedy axiom} \sl
For each $r$, ${\mathcal A}_r\co{\mathcal M}_0 \longrightarrow \{0,1\}$ satisfies
the greedy axiom (2) in definition \ref{def: local arbiter}.
In other words, given $(M,{\gamma}),(M',{\gamma}')\in {\mathcal M}_0$ with $(M,{\gamma})\subset (M',{\gamma}')$, if
${\mathcal A}_r(M,{\gamma})=1$ then ${\mathcal A}_r(M',{\gamma}')=1$.
\end{lemma}

The proof follows from the following proposition.
Recall from the beginning of section \ref{sec: Local arbiters}
that each element $(M,{\gamma})\in {\mathcal M}$ comes equipped with an
embedding $(M,{\gamma})\subset (D^4,S^3)$.

\begin{proposition} \label{proposition}  \sl \nl
(1) If $I<J$ then $(A_I,
{\alpha}_I)\not\subset (A_J, {\alpha}_J)$.

(2)  For any $I, J,$ $(A_I, {\alpha}_I)\not\subset (B_J, {\beta}_J)$.

(3)
For any multiindices $I,J$ whose first entries equal $1$, there does not exist an embedding $(B_I,
{\beta}_I)\subset (A_J, {\alpha}_J)$.
\end{proposition}

We postpone the proof of proposition \ref{proposition} to complete the
proofs of lemma \ref{greedy axiom} and theorem \ref{uncountable}. One needs to check that for
any two elements $(M_1,{\gamma}_1), (M_2,{\gamma}_2)\in {\mathcal M}_0$ if ${\mathcal A}_r(M_1)=1$ and ${\mathcal
A}_r(M_2)=0$ then
$(M_1,{\gamma}_1)\not\subset (M_2,{\gamma}_2)$. Consider the list of all pairs $M_1, M_2$ where
${\mathcal A}_r(M_1)=1$ and ${\mathcal A}_r(M_2)=0$:

(a) $M_1=A_I$ where $I<r$, $M_2=A_J$ where $J>r$,

(b) $M_1=A_I$ where $I<r$, $M_2=B_J$ where $J<r$,

(c) $M_1=B_I$ where $I>r$, $M_2=A_J$ where $J>r$,

(d) $M_1=B_I$ where $I>r$, $M_2=B_J$ where $J<r$.

The cases (a), (b), (c) follow from the statements (1), (2), (3) of proposition \ref{proposition}
respectively. To eliminate (d), suppose
$B_I\subset B_J$ where $J<I$. Considering their complements, there
is an inclusion  $A_J\subset A_I$, a contradiction with
proposition \ref{proposition} (1). This completes the proof of lemma \ref{greedy axiom}. \qed

\begin{figure}[ht]
\medskip
\includegraphics[height=6cm]{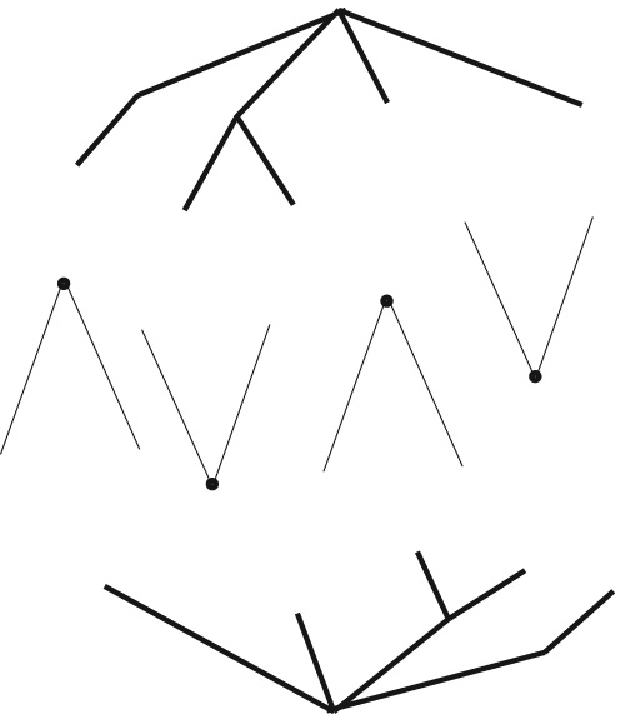}
{\scriptsize
    \put(-70,177){$A=2-$handle}
    \put(-145,108){$B_I$}
    \put(-66,105){$A_J$}
    \put(-106,44){$A_I$}
    \put(-27,70){$B_J$}
    \put(-145,64){${\mathcal A}_r\!=\!0$}
    \put(-112,91){${\mathcal A}_r\!=\!1$}
    \put(-73,-14){$B=$collar}}
\medskip
     \caption{The vertices of the graph encode codimension zero
submanifolds of $D^4$ (elements of ${\mathcal M}$), and the edges
correspond to inclusions. This figure illustrates the case
$I<r$, $J>r$: then ${\mathcal A}_r(A_I)={\mathcal A}_r(B_J)=1,$ and
${\mathcal A}_r(B_I)={\mathcal A}_r(A_J)=0$. These requirements
force ${\mathcal A}_r(M)=1$ for any $M$ in the ``cones'' over $A_I,
B_J$ and ${\mathcal A}_r(M)=0$ for any $M$ in the ``cones'' under
$B_I, A_J$. The upward-facing cones are disjoint from the
down-facing cones by lemma \ref{greedy axiom}.}
\smallskip
\label{fig:partial ordering}
\end{figure}

\subsection{Proof of theorem \ref{uncountable}.} \label{sec:proof}
Given a ray $r$, by lemma \ref{greedy axiom} the partial
arbiter ${\mathcal A}_r$ is consistently defined on the collection ${\mathcal M}_0$
of model decompositions. Given $M_1, M_2\in {\mathcal M}_0$  with ${\mathcal
A}_r(M_1)=0$, ${\mathcal A}_r(M_2)=1$, axiom (2) in definition \ref{def: local arbiter} forces ${\mathcal
A}_r(M)=0$ for any $M$ isotopic to a submanifold of $M_1$ and ${\mathcal A}_r(M)=1$ for any
$M$ isotopic to $M'$ where $M_2\subset M'$. Graphically, ${\mathcal A}_r=0$ in the ``cone''
below $M_0$ and ${\mathcal A}_r=1$ in the ``cone'' above $M_1$,
figure \ref{fig:partial ordering}. The upward-facing cones are disjoint from the down-facing
cones by lemma \ref{greedy axiom}. This defines a consistent extension of ${\mathcal
A}_r$ to the union of all such cones. Duality axiom (3) in definition \ref{def: local arbiter}
for this extension follows from the fact that the upward-facing cone
over, say, $A_I$ consists precisely of the complements in $D^4$ of
the elements of the down-facing cone under $B_I$.
To extend ${\mathcal A}_r$ to ${\mathcal M}$, consider any local arbiter on $D^4$, for example a homological
arbiter ${\mathcal A}_F$ for some fixed field $F$. Set ${\mathcal A}_r(M)={\mathcal A}_F(M)$ for any $M$
in the complement of the union of the ``cones'' discussed above. It is straightforward to check that ${\mathcal A}_r$
satisfies the axioms of a local arbiter.

Note that if $r\neq r'$ then ${\mathcal A}_r\neq {\mathcal A}_{r'}$:
assume $r<r'$, then their difference is detected by a codimension zero submanifold $A_I$ where
$I$ is a multiindex with $r<I<r'$.  This completes the proof of theorem \ref{uncountable}, assuming
proposition \ref{proposition}. \qed

{\bf Remarks.}\\
{\bf 1.} Because we use a definite procedure to extend ${\mathcal A}_r$ beyond the ``cones''
on which it is determined by $r$, the proof of theorem \ref{uncountable} does not depend on the axiom of choice.
Nor does any other assertion in this paper. To be complete, our argument actually has shown
that the cardinality of local arbiters on $D^4$ is precisely the continuum.

{\bf 2.} The argument extending partial arbiters ${\mathcal A}_r$ to arbiters defined
on all codimension zero submanifolds of $D^4$ in the proof of theorem \ref{uncountable} may also be used to define ``hybrids''
of homological arbiters of different characteristics. For example, suppose $(M, {\gamma})$ is an element of
${\mathcal M}$ such that ${\gamma}=0\in H_1(M;{\mathbb Z}/p)$ but ${\gamma}\neq 0\in H_1(M;{\mathbb Z}/q)$
for some primes $p$ and $q$. By Alexander duality with ${\mathbb Z}/p$ coefficients, two copies
of $M$ cannot be embedded disjointly into $D^4$ so that the distinguished curves form the Hopf link on the boundary.
Then a partial arbiter ${\mathcal A}$ may be defined by setting it equal to
${\mathcal A}_{{\mathbb Z}/p}$ in the cone over $M$ and also in the cone under $D^4\smallsetminus M$. Extend it by
${\mathcal A}_{{\mathbb Z}/q}$ to the rest of ${\mathcal M}$. The arbiters constructed in the proof of theorem
\ref{uncountable} assume value $1$ on elements $(M, {\gamma})\in {\mathcal M}$ such that
${\gamma}\neq 0\in H_1(M;R)$ for {\bf any} coefficients $R$, therefore they
are different from any such hybrid of homological arbiters.

\smallskip

{\bf Proof of part (1) of proposition \ref{proposition}.} Before giving the
general (inductive) proof consider several special cases which
clarify the inductive step. (Refer to figure \ref{fig:model decompositions}
and the beginning of section \ref{sec: Model decompositions} to recall the definition
of the submanifolds $A_I, B_I$.)

(i) {\sl $(A_0, {\alpha}_0)$ does not embed in $(A_1, {\alpha}_1)$.}

The proof is immediate since ${\alpha}_0=0\in H_1(A_0)$ and
${\alpha}_1\neq 0\in H_1(A_1)$.

(ii) {\sl $(A_{10}, {\alpha}_{10})\not\subset (A_{11}, {\alpha}_{11})$.}

Suppose such an embedding exists. Then consider the complements:
$(B_{11},{\beta}_{11})\subset (B_{10},{\beta}_{10})$.
The submanifold $B_{11}$ is (a thickening of) a $2$-stage {\em grope}, see \cite{FL} (the term used in that reference is
a {\em half-grope}), a geometric counterpart of the commutators in the lower central series ${\pi}_n$ of a group $\pi$. This series
is defined inductively by ${\pi}_n=[{\pi},{\pi}_{n-1}]$, ${\pi}_1={\pi}$. The attaching curve ${\beta}_{11}$
of the $2$-stage grope $B_{11}$ is in the
third term $({\pi}_1 (B_{11}))_3$.
This is a contradiction since ${\beta}_{11}$ maps to ${\beta}_{10}$, however
$${\beta}_{11}=1\in{\pi}_1(B_{11})/({\pi}_1(B_{11}))_3 \hspace{.25cm} {\rm
but} \hspace{.25cm} {\beta}_{10}\neq
1\in{\pi}_1(B_{10})/({\pi}_1(B_{10}))_3.$$

(iii) {\sl $(A_{00},
{\alpha}_{00})\not \subset (A_{01}, {\alpha}_{01})$.}

\begin{figure}[ht]
\vspace{.5cm}
\includegraphics[height=4cm]{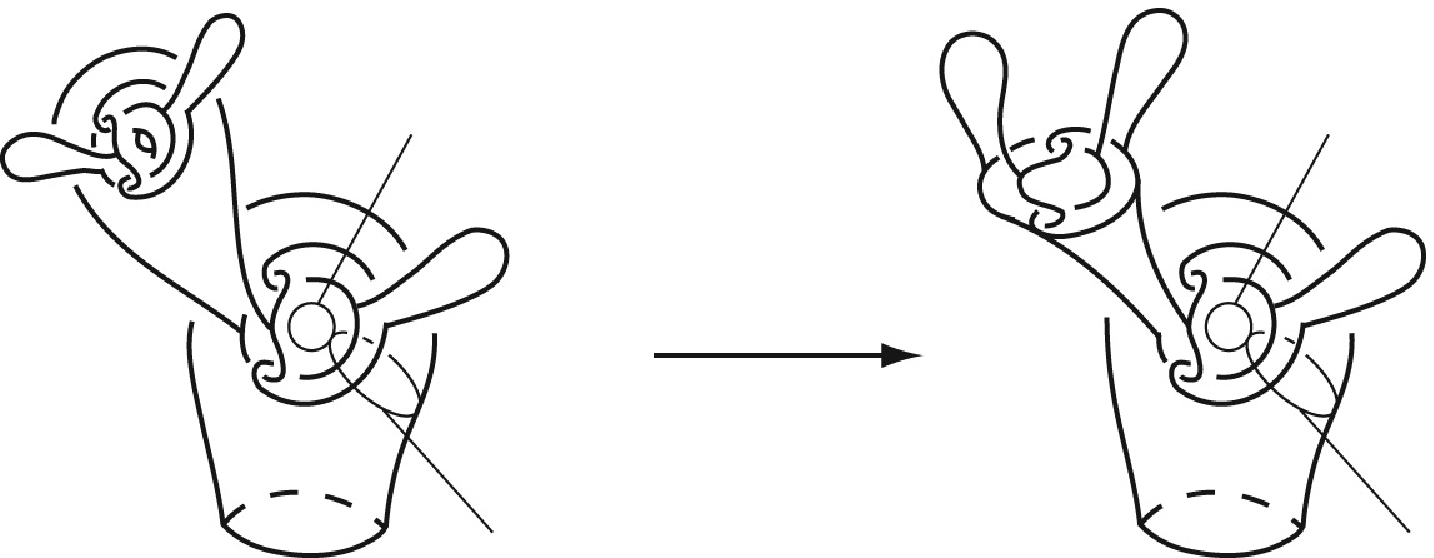}
{\footnotesize
    \put(-283,42){$A_{00}$}
    \put(-137,48){$i$}
    \put(-5,35){$A_{01}$}
    \put(-270,15){${\Sigma}'$}
    \put(-80,15){${\Sigma}''$}
    \put(-250,-9){${\alpha}_{00}$}
    \put(-55,-9){${\alpha}_{01}$}
    \put(-210,90){${\gamma}_1$}
    \put(-192,-2){${\gamma}_2$}
    \put(-20,90){${\delta}_1$}
    \put(-2,-2){${\delta}_2$}}
\vspace{.35cm} \caption{}
\label{fig:non-embedding1}
\end{figure}

\begin{figure}[ht]
\vspace{.5cm}
\includegraphics[height=4cm]{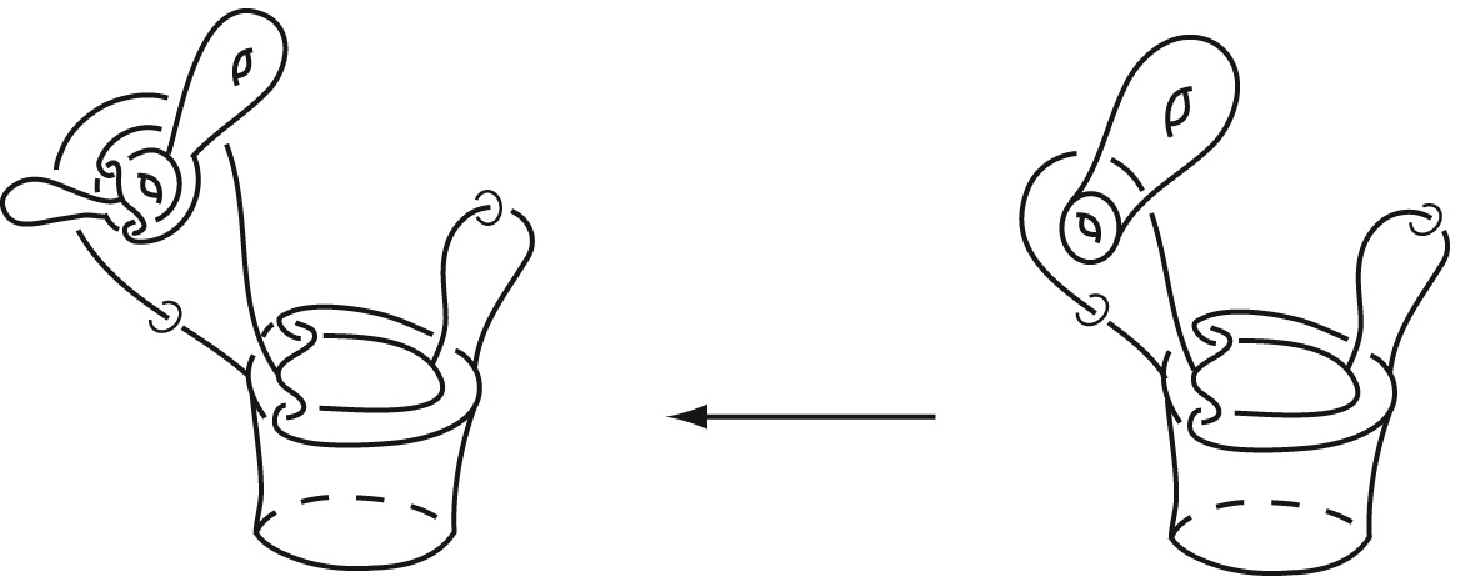}
{\footnotesize
    \put(-280,17){$B_{00}$}
    \put(-80,17){$B_{01}$}
    \put(-129,39){$\overline i$}
    \put(-250,-9){${\beta}_{00}$}
    \put(-55,-9){${\beta}_{01}$}
    \put(-265,40){${\gamma}_1$}
    \put(-196,83){${\gamma}_2$}
    \put(-80,40){${\delta}_1$}
    \put(-6,79){${\delta}_2$}
    \put(-239,68){$T'$}
    \put(-182,58){$S'$}
    \put(-102,73){$T''$}
    \put(2,58){$S''$}}
    \vspace{.35cm} \caption{}
    \label{fig:non-embedding2}
\end{figure}

The argument in this example will serve as an ingredient in the proof in the general case further below.
Suppose the embedding $i$  exists as in figure \ref{fig:non-embedding1}. Considering the complements one gets
$\overline i\co (B_{01},{\beta}_{01})\longrightarrow
(B_{00},{\beta}_{00})$, figure \ref{fig:non-embedding2}. The argument will involve two steps,
the first one using the embedding $i$ and second step using $\overline i$. The connection between the two steps
is provided by the Alexander duality.

{\em Step 1.}
Denote by ${\Sigma}', {\Sigma}''$
the bottom surface stages of $A_{00}, A_{01}$. Since
${\alpha}_{00}$ maps to ${\alpha}_{01}$, $i$ ``restricts'' to a
degree $1$ map $({\Sigma}',{\alpha}_{00})\longrightarrow ({\Sigma}'',{\alpha}_{01})$. (More
precisely, consider the restriction $i|_{{\Sigma}'}$. $A_{01}$ is
homotopy equivalent to the wedge of ${\Sigma}''$ and three
$2-$spheres. The composition of $i|_{{\Sigma}'}$ with the
retraction onto ${\Sigma}''$ has degree $1$.)

Consider the classes in $H_2(B_{00})$, $H_2(B_{01})$ represented by
the tori $T'$, $T''$ and spheres $S'$, $S''$, figure \ref{fig:non-embedding2}. (These are
closed surfaces since they are capped off by the null-homotopies for
the Bing-double curves.) The punctured tori ${\Sigma}'\subset
A_{00}, {\Sigma}''\subset A_{01}$ have preferred symplectic bases of
curves: ${\gamma}_1, {\gamma}_2$ on ${\Sigma}'$, and ${\delta}_1, {\delta}_2$ on
${\Sigma}''$. The curves ${\gamma}_1, {\gamma}_2$ are Alexander dual
respectively to $T', S'$, and similarly ${\delta}_1, {\delta}_2$
are dual to $T'',S''$. Consider the matrix $M=\left(
\begin{smallmatrix} a & b\\ c & d \end{smallmatrix} \right)$ representing
the induced map $i_*\co H_1({\Sigma}')\longrightarrow
H_1({\Sigma}'')$ in these bases. Since the map
$({\Sigma}',\partial)\longrightarrow ({\Sigma}'',\partial)$ has
degree $1$, $M$ is an element of $SL_2({\mathbb Z})$.

{\em Step 2.} The matrix for the map $\overline i_*\co H_2(B_{01})\longrightarrow H_2(B_{00})$,
restricted to the $2$-dimensional subspaces spanned by $[S''],
[T'']$ and $[S'], [T']$, is the transpose of $M$, by Alexander
duality. Since $[S'']$ is a spherical class, the component of
$\overline i_*([S''])$ along $[T']$ is trivial, so $b=0$ and $M$
must be of the form $\left(
\begin{smallmatrix} 1 & 0\\ c & 1 \end{smallmatrix} \right)$.
Therefore the class $[T'']$ maps to $[T']$ plus a multiple $c$ of $[S']$,
but in particular one gets a degree $1$ map $T''\longrightarrow T'$.
This leads to a
contradiction when one considers the elements represented by $T'$, $T''$ in
the {\em Dwyer's filtration} (see \cite{Dwyer}) of the corresponding spaces.
A geometric description of the $n$th term of the Dwyer's filtration ${\phi}_n(X)
\subset H_2(X;{\mathbb Z})$ of a space $X$ is the subgroup consisting of those
elements in $H_2(X;{\mathbb Z})$ which may be represented by maps of closed gropes
of class $n$ into $X$.
Since $[T'']\in {\phi}_3(B_{01})$ but $[T']\not\in
{\phi}_3(B_{00})$, $[T'']$ cannot map to $[T']$, a contradiction.

\medskip

{\bf Proof of part (1) of proposition \ref{proposition} in the general case.} Suppose an
embedding $i\co (A_I, {\alpha}_I)$ $\longrightarrow (A_J,
{\alpha}_J)$ exists; let $\overline i\co
(B_J,{\beta}_J)\longrightarrow (B_I,{\beta}_I)$ be the corresponding
inclusion of the complements. Since $I<J$, $I=(i_1,\ldots, i_n, 0,
\ldots)$ and $J=(i_1, \ldots, i_n, 1, \ldots)$ for some $n$.

First consider the case $i_1=\ldots=i_n=0$, the proof for arbitrary $i_1,\ldots,i_n$ will be analogous.
Since all $i_j$ are assumed to be zero, $j=1\dots,n$, each of the first $n$ steps in the construction
of $A_I, A_J$ replaces an existing $2$-handle with (a thickening of) the genus $1$ surface, and then attaches
a pair of $2$-handles to the Bing double of an essential curve on this surface. Denote the
successive genus $1$ surface stages of $A_I$ by ${\Sigma}'_1,\ldots, {\Sigma}'_n$, and similarly
denote the genus $1$ surface stages of $A_J$ by ${\Sigma}''_1,\ldots, {\Sigma}''_n$.
The construction on the $B$-side is dual (in the sense of proposition \ref{surface complement}):
each of the first $n$ steps attaches a pair of $2$-handles
to a Bing double, and then replaces one of these $2$-handles with a punctured torus. Denote the
successive tori in $B_I$ by $T'_1,\ldots, T'_n$ and the corresponding tori in $B_J$ by  $T''_1,\ldots, T''_n$.
(In the special case $I=(00), J=(01)$ considered in (iii) above, $n$ was equal to $1$.)

{\em The base of the induction}.
The embedding $i$ induces a degree $1$ map $({\Sigma}'_1,{\alpha}_I)
\longrightarrow ({\Sigma}''_1,{\alpha}_J)$.

{\em The inductive step}.
The inductive hypothesis is that the embedding $i$ induces a degree $1$ map ${\Sigma}'_k
\longrightarrow {\Sigma}''_k$ for some $1\leq k<n$. The argument in the special case (iii)
above proves that then the complementary embedding $\overline i$ induces a degree $1$ map $T'_k
\longrightarrow T''_k$. Since the construction of the first $n$ stages of the submanifold
$A_I$ is precisely dual to the construction of $B_I$ (and similarly for $A_J, B_J$, compare figures
\ref{fig:non-embedding1}, \ref{fig:non-embedding2}), exactly the same
argument now propagates the conclusion to the $A$-side: if $\overline i$ induces a degree $1$ map
$T''_k \longrightarrow T'_k$, then the embedding $i$ induces a degree $1$ map ${\Sigma}'_{k+1}
\longrightarrow {\Sigma}''_k$. This concludes the proof of the inductive step.

A contradiction is reached after $n-1$ applications of the inductive
step: the outcome is a degree $1$ map ${\Sigma}'_n \longrightarrow {\Sigma}''_n$;
run the argument one more time to get a degree $1$ map
$T''_n\longrightarrow T'_n$. However since the $n+1$-st coefficient of $J$ is $1$, it follows
from the construction of $B_J$ that $[T''_n]\in {\phi}_3(B_J)$. On the other hand
the $n+1$-st coefficient of $I$ is $0$, so $[T'_n]\not\in {\phi}_3(B_I)$. This
contradicts the fact that $\overline i$ induces a degree $1$ map $T''_n\longrightarrow T'_n$, concluding the
proof when $i_1=\ldots=i_n=0$.

A slight modification of the argument is necessary in the general case when $i_1=0$ but
$(i_1,\ldots,i_n)$ contains subsequences with one or more consecutive $1$'s.
Then the corresponding tori on the $B$-side are replaced by
gropes of height $h$ equal to ($1+$ the number of consecutive $1$'s in a subsequence),
all of whose surface stages have genus $1$.  (Correspondingly, the $2$-handles on the $A$-side are attached
to higher iterated Bing doubles.) The inductive step discussed above applies without any changes on the $A$-side.

On the $B$-side a variation is necessary: Let $f$ be a map $f\co G\longrightarrow G$ where $G$ is a half-grope of
height $n$ with all surface stages of genus $1$. The second homology class of the bottom stage surface
of $G$ reads off the relation $[{\gamma}_1,[{\gamma}_2,\ldots,[{\gamma}_h,{\gamma}_{h+1}]\ldots]$, where
${\gamma}_i$ are the curves generating ${\pi}_1(G)$, corresponding to the ``tips'' of the grope.
Note that ${\gamma}_h,{\gamma}_{h+1}$ generate ${\pi}_1$ of the top stage surface of $G$.
Suppose $f_*$ carries the class of the bottom stage surface of $G$ in $H_2(G)$ to itself. Then the map
induced by $f$ on the $2$-dimensional subspace of $H_1(G)$ spanned by ${\gamma}_h,{\gamma}_{h+1}$ is in
$SL(2,{\mathbb Z})$. The rest of the argument in the inductive step goes through.

Finally, if $i_1=1$ then the first step of the induction is applied on the
$B$-side. This concludes the proof of \ref{proposition} (1).

\medskip

{\bf The proof of part (2) of  proposition \ref{proposition}} involves a generalization of the Milnor group
from the context of link homotopy (a link in $S^3$ bounding disjoint maps of disks in $D^4$)
to a class of more general submanifolds of the $4$-ball.
We outline this construction below and refer the reader to \cite{Krushkal} for a detailed
exposition.

Suppose $(A_I,{\alpha}_I)\subset (B_J, {\beta}_J)$.
Then $A_I$ embeds in the complement of $A_J$, where the attaching curves
${\alpha}_I$, ${\alpha}_J$ form the Hopf link in $S^3=\partial D^4$.
As a warmup, consider the special case where both $I, J$ are dyadic multiindices with no zeros, for example
$A_{1}$ and $A_{11}$ are shown in figure \ref{fig:model decompositions}. For such multiindices $I, J$, the
handlebodies $A_I, A_J$ consist of $2$-handles attached to iterated Bing doubles of their respective
attaching curves ${\alpha}_I, {\alpha}_J$ which form the Hopf link in $\partial D^4$. It is straightforward
to prove that any link $L$ obtained by iterated Bing doubling of the Hopf link is a homotopically essential link \cite{Milnor},
by computing that it has a non-trivial Milnor's ${\bar\mu}$-invariant with non-repeating coefficients. Therefore the components of
$L$ do not even bound disjoint maps of disks into $D^4$, contradicting the assumption
$(A_I,{\alpha}_I)\subset (B_J, {\beta}_J)$.

\begin{figure}[ht]
\centering
\bigskip
\includegraphics[width=3.6cm]{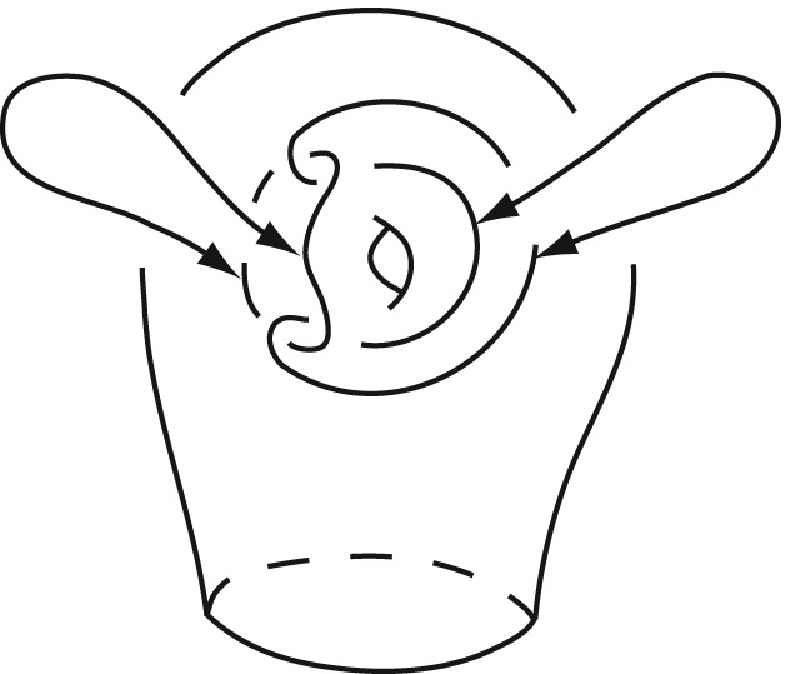}
 \hspace{1.5cm} \includegraphics[width=5.6cm]{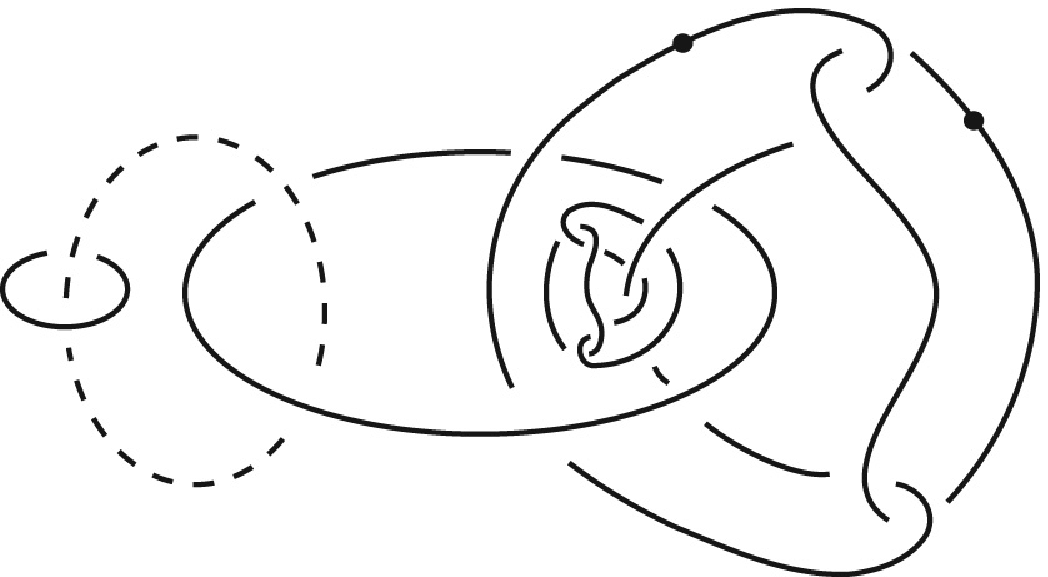}
{\small
    \put(-309,36){$A_0$}
    \put(-291,-3){${\alpha}_0$}
    \put(-175,42){${\alpha}_0$}}
{\scriptsize
    \put(-99,13){$0$}
    \put(-82.5,42){$0$}
    \put(-53,42){$0$}}
    \smallskip
\caption{The submanifold $A_0$ and its Kirby diagram. (The handle diagram is
drawn in the solid torus which is the complement of the dashed circle in the $3$-sphere.)}
\label{fig:A0}
\end{figure}

\begin{figure}[ht] \label{model figure1}
\centering
\bigskip
\includegraphics[width=5.75cm]{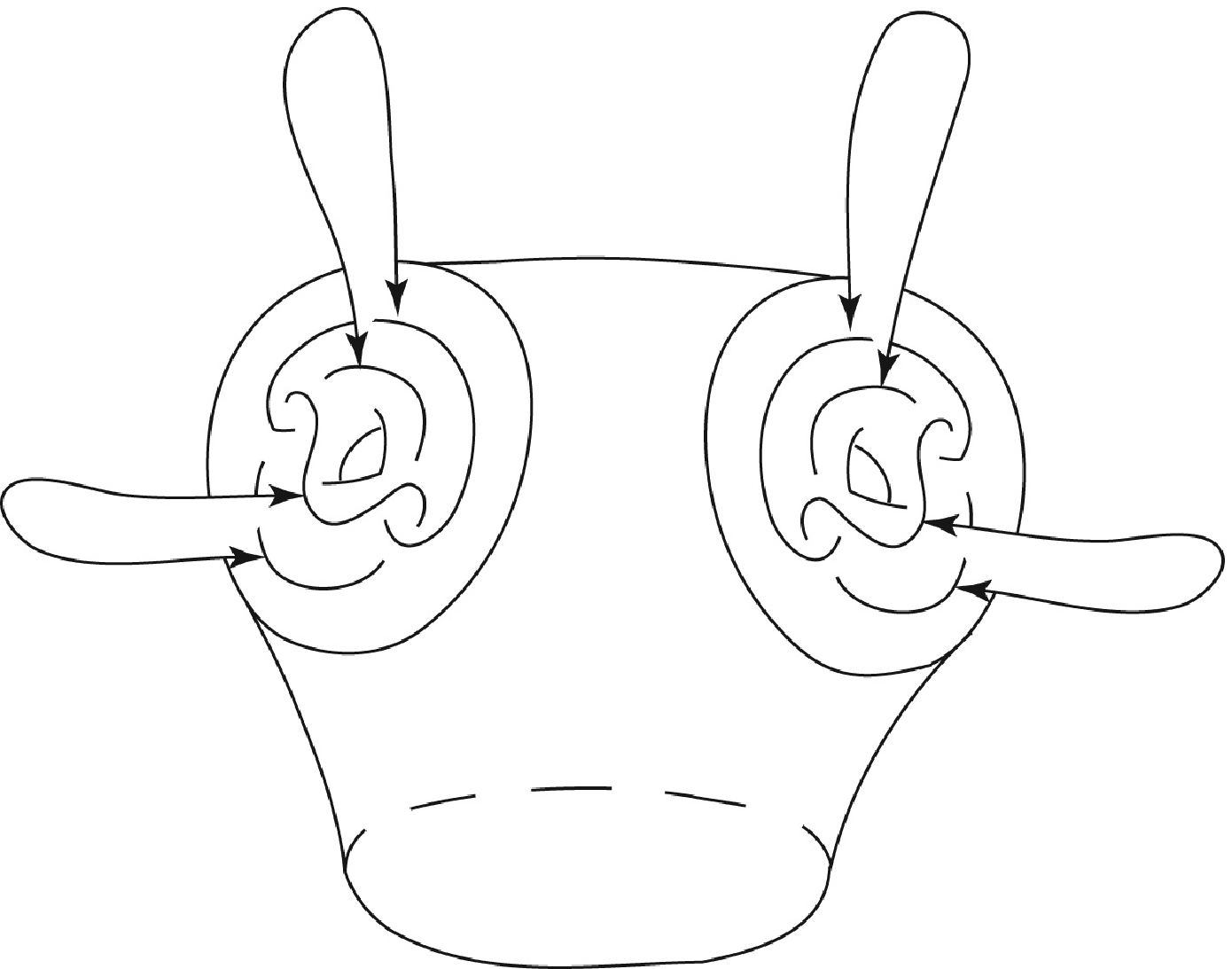}\hspace{2cm} \includegraphics[bb=0 -55 250 100, width=5.5cm]{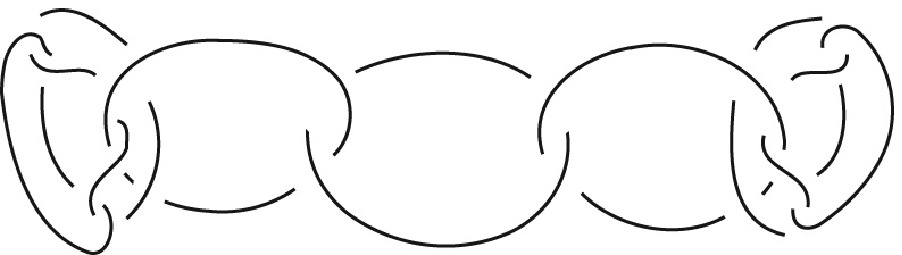}
{\small
    \put(-349,9){${\gamma}$}
    \put(-387,73){$h_1$}
    \put(-256,20){$C$}
    \put(-360,122){$h_2$}
    \put(-282,120){$h_3$}
    \put(-227,65){$h_4$}}
{\scriptsize
    \put(-164,58){$0$}
    \put(-141,80){$0$}
    \put(-18,79){$0$}
    \put(7,63){$0$}}
    {\small
    \put(-75,25){${\gamma}$}}
    \put(-118,73){\circle*{3}}
    \put(-39,72){\circle*{3}}
    \vspace{.3cm}
\caption{The handlebody $C$ (a Bing cell of height $1$): a schematic
picture and a Kirby diagram}
\label{cell figure}
\end{figure}

Now consider a more interesting example $A_0$ (figure \ref{fig:A0}), illustrating the argument in the general case.
Consider a pair of pants $P$ with boundary components ${\gamma}, {\alpha}_1, {\alpha}_2$, and let
$C$ denote $P\times D^2\; \cup$ four zero-framed $2$-handles $h_1\ldots,h_4$ attached to the Bing doubles of the curves ${\alpha}_1, {\alpha}_2$,
figure \ref{cell figure}.

Observe that there is a map $C\longrightarrow A_0$ taking $\gamma$ to ${\alpha}_0$, such that the $2$-handle $h_1$ is disjoint from $h_2$,
and similarly $h_3$ is disjoint from $h_4$. On the level of spines, the proof amounts to cutting the punctured torus into
a pair of pants; the separation of parallel copies of the $2$-handles results in $h_1-h_3$ and $h_2-h_4$ double points.

Therefore to prove that the components of the Hopf link do not bound disjoint embeddings of two copies of $A_0$, it suffices to show that
its components do not bound disjoint maps of two copies $C', C''$ of $C$, subject to the disjointness condition on the $2$-handles discussed above.
Suppose the opposite is true. Consider meridians $m'_i$ to
the $2$-handles $\{h'_i\}$ of $C'$ and meridians $m''_i$ to the $2$-handles $\{ h''_i\}$ of $C''$ in $D^4$.
(By a meridian we mean a based loop in the complement representing the same homology class as a fiber
of the normal circle bundle over a $2$-handle.) Consider the quotient $\overline M {\pi}$
of ${\pi}:={\pi}_1(D^4\smallsetminus (C'\cup C''))$ by the relations $[(m'_i)^x,( m'_j)^y]$ and $[(m''_i)^x, (m''_j)^y]$, where
$\{i,j\}\neq \{1,2\}, \{3,4\}$ and the conjugating elements $x,y$ range over all elements of ${\pi}$.
These relations correspond to the Clifford tori linking the double points between
the various $2$-handles in $D^4$, and the fact that the pairs $\{1,2\}, \{3,4\}$ are excluded reflects the fact that these
handles are required to be disjoint. (This is a generalization of the Milnor group $M{\pi}$ (see \cite{Milnor}) which in this context is defined
as the quotient of ${\pi}$ by the relations  $[(m'_i)^x,( m'_i)^y], [(m''_i)^x,( m''_i)^y], i=1,\ldots, 4$, which correspond to
{\em self}-intersections of the handles. These relations are included in the defining relations of $\overline M{\pi}$.)

The classical Milnor group $M{\pi}$ is nilpotent, and since $\overline M{\pi}$ is a quotient of $M{\pi}$, it is nilpotent
as well. It follows that $\overline M{\pi}$ is generated by the chosen meridians, and its relations may be read off from the
generators of $H_2(D^4\smallsetminus (C'\cup C''); {\mathbb Z})$, which is Alexander dual to the relative
group $H_1(C'\cup C'', \partial(C'\cup C'')\cap \partial D^4)$. There are two types of generators of this second homology
group: the Clifford tori that give rise to the commutation relations used in the definition of $\overline M{\pi}$, and
the generators dual to the relative first homology group of the pairs of pants defining $C', C''$. The latter relations
are of the form $[m'_1,m'_2]=[m'_3,m'_4], [m''_1,m''_2]=[m''_3,m''_4]$.

Let $m', m''$ denote the respective meridians to the components of the Hopf link in $\partial D^4$,
formed by the attaching curves ${\gamma}', {\gamma}''$. The relation $[m', m'']=1$ holds in ${\pi}_1(S^3\smallsetminus
({\gamma}'\cup {\gamma}''))$, therefore it also holds in ${\pi}_1(D^4\smallsetminus (C'\cup C''))$. On the other hand,
in terms of the generators of $\overline M{\pi}$, $m'=[m'_1,m'_2]=[m'_3,m'_4], m''=[m''_1,m''_2]=[m''_3,m''_4]$.
One then shows using the Magnus expansion that, based on the relations in $\overline M{\pi}$ discussed above,
the commutator $[m', m'']$ is in fact non-trivial in $\overline M{\pi}$: phrased differently, there are
``not enough'' relations in ${\pi}_1(D^4\smallsetminus (C'\cup C''))$ to imply the relation $[m', m'']=1$
which holds in ${\pi}_1(S^3\smallsetminus ({\gamma}'\cup {\gamma}''))$.
This contradiction shows that
the components of the Hopf link do not bound in $D^4$ disjoint copies of the manifold $A_0$.

The manifolds $A_0$, $A_1$ considered above are the building blocks in the definition of model decompositions
$D^4=A_I\cup B_I$, see section \ref{sec: Model decompositions}. There is a similar inductive construction of {\em Bing cells}
based on the manifold $C$ above,
where the inductive step replaces a $2$-handle with a copy of $C$. {\em Moreover, the $2$-handles in the more general definition
below may be attached to iterated} Bing doubles, since such links arise in the model submanifolds $A_I$ (see for example
$A_{01}$ in figure \ref{fig:model decompositions}.)

\begin{definition} \label{model cell}
A {\em model Bing cell of height $1$} is a smooth $4$-manifold $C$
with boundary and with a specified curve ${\gamma}\subset
\partial C$, defined as follows.
Consider a pair of pants $P$ with boundary components ${\gamma}, {\alpha}_1, {\alpha}_2$, and let
$L_i$ be a (possibly iterated) Bing double of the core of the solid torus ${\alpha}_i\times D^2$.
Then $C$ is obtained from
$P\times D^2$ by attaching zero-framed $2$-handles along the
components of $L_1\cup L_2$.

A {\em model Bing cell $C$ of height $h$} is obtained from a model
Bing cell of height $h-1$ by replacing some of its $2$-handles with Bing cells
of height one.
\end{definition}

{\em Bing cells in a $4$-manifold} $M$ are defined as maps of model
Bing cells into $M$, subject to certain crucial disjointness
requirements. For each link $L_i$ in the inductive definition above,
surfaces and more generally Bing cells attached to different components of
$L_i$  are required to be disjoint in $M$.
Generalizing the arguments above, one shows  (see \cite{Krushkal}) that

(1) \parbox[t]{14.25cm}{If the components of a link $L\subset S^3=\partial D^4$ bound
disjoint Bing cells in $D^4$ then $L$ is homotopically trivial, and}

(2) For each index $I$ the attaching curve ${\alpha}_I$ bounds a Bing cell in $A_I$.

Therefore the components of the Hopf link do not bound disjoint embeddings of $A_I, A_J$, i.e. $(A_I,{\alpha}_I)\not
\subset (B_J, {\beta}_J)$. This concludes the proof of proposition \ref{proposition} (2).

{\em The proof of part (3) of proposition} \ref{proposition} is immediate since
${\beta}_I=0\in H_1(B_I)$ and ${\alpha}_J\neq 0\in H_1(A_J)$.
This completes the proof of
proposition \ref{proposition} and of theorem \ref{uncountable}. \qed

\subsection{Arbiters and embeddings.} \label{sec: embeddings}

In this section we note the important role played by the embeddings in definition
\ref{def: local arbiter} of arbiters on $D^4$.
According to \cite{K0} for any given $g$ there exist decompositions (in the sense of axiom (3) in definition
\ref{def: local arbiter}) $D^4=A_g\cup B_g$, where the attaching curves ${\alpha}, {\beta}$ form the Hopf link
in $\partial D^4$, such that the conditions (i), (ii) are satisfied:

(i) Let $S$ denote the genus $g$ surface with one boundary
component, ${\gamma}=\partial S$. Denote by $S_0$ its
untwisted $4-$dimensional thickening, $S_0=S\times D^2$, and
set ${\gamma}_0={\gamma}\times\{0\}$. Then there exists a
proper embedding $(A_g, {\alpha})\subset (S_0,
{\gamma}_0)$.

(ii) $B_g$ embeds in a collar on its attaching curve. More precisely, there
exists a proper embedding $(B_g, {\beta})\subset
(S^1\times D^2\times [0,1],\, S^1\times \{0\}\times\{0\})$

Suppose ${\mathcal A}$ is an invariant as in definition \ref{def: local arbiter}, where the embedding in the axiom (2) in \ref{def: local arbiter}
is assumed to be an arbitrary embedding of abstract pairs ($M,{\gamma}), (M', {\gamma}')$, rather than just the inclusion of submanifolds of $D^4$. Since the value of ${\mathcal A}$ on the collar is necessarily $0$, (ii) implies
that for each $g$, ${\mathcal A}(B_g)=0$, therefore ${\mathcal A}(A_g)=1$. (Note that the embeddings in (i), (ii) are indeed not isotopic to
the original embeddings into $D^4$.)
Now suppose $(M,{\gamma})$ is a submanifold of $(D^4, S^3)$ such that the attaching curve $\gamma$
of the $4$-manifold $M$ is trivial in $H_1(M;{\mathbb Z})$. Then $\gamma$ bounds an embedded surface $S$ in $M$ with the trivial
normal bundle. It follows that $(A_g, {\alpha})\subset (M,{\gamma})$, where $g$ is the genus of $S$. Again by axiom (2) which allows
arbitrary embeddings,
${\mathcal A}(M,{\gamma})=1$. That is, a function ${\mathcal A}$ which is insensitive to the embedding of submanifolds into $D^4$
is determined by homology, at least on those submanifolds whose attaching curve is trivial in the first homology with integer coefficients, and
on their complements in $D^4$.
The proof of theorem \ref{uncountable} is valid only when the embeddings are inclusions of submanifolds of $D^4$, indeed
the arbiters constructed in its proof assume value $1$ on manifolds where the attaching circle does not vanish in homology with any
coefficients.

\subsection{Submanifolds of ${\mathbb C}P^2$.} \label{sec: CP2}

This section outlines a construction of submanifolds of ${\mathbb C}P^2$ which may be viewed as
an ``absolute'' analogue of the relative model decompositions of $D^4$ defined in section \ref{sec: Model decompositions}.

Consider the usual handle structure on ${\mathbb C}P^2$ with one handle of each index $0,2,4$:
${\mathbb C}P^2=0H\cup2H\cup 4H$. Consider decompositions ${\mathbb C}P^2=\widehat A\cup \widehat B$, where
$\widehat A\cap \widehat B=\partial \widehat A\cap \partial \widehat B$. The decompositions are constructed of the form
$\widehat A=0H\cup A$, $\widehat B=4H\cup B$, and $2H=A\cup B$ is a decomposition of the $4$-ball
$2H=D^2\times D^2$ in the sense discussed in section \ref{sec: Local arbiters}: $A\cup B$
extends the standard genus $1$ decomposition of the $3$-sphere $D^2\times S^1\cup S^1\times D^2$.

Applying this construction to model decompositions $D^4=A_I\cup B_I$ of the $4$-ball in section \ref{sec: Model decompositions}
(see figure \ref{fig:model decompositions}), one gets a dyadic collection of decompositions ${\mathbb C}P^2=\widehat A_I\cup \widehat B_I$.
More concretely, $\widehat A_I$ is obtained by attaching a $1$-framed $2$-handle to $A_I$ along the attaching curve ${\alpha}_I$, and
similarly its complement $\widehat B_I$ is obtained by attaching a $1$-framed $2$-handle to $B_I$ along the attaching curve ${\beta}_I$.

It seems reasonable to conjecture that there are non-embedding results analogous to proposition \ref{proposition}, which should
enable one to construct an uncountable collection of topological arbiters on ${\mathbb C}P^2$. However proving this would require an
extension of the Milnor group techniques which is outside the scope of this paper. This question remains open.

\section{Squares in stable homotopy and arbiters in higher dimensions} \label{sec: Higher dimensions}

Recall that homology with field coefficients gives rise to topological arbiters on even dimensional projective spaces,
and similarly there are local homological arbiters on $D^{2n}$ for all $n$, see lemma \ref{projective space} and the
discussion following definition \ref{def: local arbiter}.
In this section we use stable homotopy theory to show that there exist local arbiters different from the homological ones
in even dimensions $\geq 6$.

Recall that the construction of an uncountable
collection of local arbiters on $D^4$ in section \ref{sec: Local arbiters}
was based on obstructions related to the failure of the Whitney trick in dimension $4$. More specifically, this is reflected in
the fact that the middle-dimensional (index $2$) handles of certain submanifolds of $D^4$ cannot be embedded disjointly. This argument
does not have an analogue in higher dimensions, since the Whitney trick is valid then.
Our construction in dimensions $6$ and higher, discussed below, is quite different. It is homotopy-theoretic in nature, and
the non-embedding results involve handles above the middle dimension. We point out that conversely, this higher-dimensional
approach does not work in dimension $4$ since the higher homotopy groups of the circle (the attaching region of the
submanifolds in $D^4$) are trivial.

Before outlining a general construction of families of non-homological arbiters on $D^{2n}$ for all $n\geq 4$, we discuss a specific
example in dimension $8$. (The construction in $D^6$ requires an additional step, discussed in section \ref{D6} below.) Consider codimension zero submanifolds $M\subset D^8$ such that $\partial M \cap \partial D^8$ is a regular
neighborhood of a standard $S^3\subset S^7=\partial D^8$. We will call this $S^3\subset \partial M$ the attaching $3$-sphere of $M$,
and consider {\em proper} embeddings $M\hookrightarrow D^8$ taking the attaching $3$-sphere to $\partial D^8$. Our construction
of a non-homological arbiter on $D^8$ relies on the following result.

\begin{lemma} \label{Hopf lemma} \sl
There exist submanifolds $(M,S^3)\subset (D^8,S^7)$ such that

(1) \parbox[t]{14.25cm}{for any coefficient ring $R$ the attaching $3$-sphere
of $M$ is non-trivial in $H_3(M; R)$, and}

(2) \parbox[t]{14.25cm}{The components of the Hopf link $S^3\cup S^3\subset S^7$ do not bound
disjoint embeddings of two copies of $M$ into $D^8$.}
\end{lemma}

{\em Proof.} Consider a standard $3$-sphere $X=S^3\subset \partial D^8$,
let $C=S^3\times D^4\times I\subset D^8$ be a collar on its regular neighborhood,
and denote by $Y$ a copy pushed into $D^8$, $Y=S^3\times\{0\}\times\{1\}$.
Let ${\Sigma}h\co S^4\longrightarrow Y$ be the suspension of the Hopf map
$h\co S^3\longrightarrow S^2$. Consider
a general position map $f\co D^5\longrightarrow D^8$ whose restriction $f|_{\partial D^5}$ is ${\Sigma}h$. Define $M$ to be the collar
$C$ union a regular neighborhood of the
image of $f$ in $D^8$.

To check (1), first note that the attaching $3$-sphere is non-trivial in homology of the abstract CW-complex
$Z:=S^3\cup_{{\Sigma}h} D^5$. The singular set $S$ of a general position map $f$ is $2$-dimensional. It follows
that the attaching $3$-sphere $X$ is non-trivial in $H_3(M)$, indeed since the third homology of the singular set
is trivial it follows from the Mayer-Vietoris sequence that $X\neq 0\in H_3(Z\cup_S cone(S))$.

Suppose the components $S^3_i$ of the Hopf link bound disjoint copies $M_i, i=1,2$ of $M$ in $D^8$.
The attaching $4$-sphere of the $5$-cell in the definition of $M$ is mapped by ${\Sigma}h$ to the copy $Y$ of $S^3$
pushed into $D^8$. Using the collar structure, it also may be viewed as being mapped by ${\Sigma}h$ to the
actual attaching sphere $X\subset \partial D^8$. Consider the resulting maps $S^4_i\longrightarrow X$.

It is convenient to consider $S^7=\partial D^8$ as ${\mathbb R}^7$ union the point at infinity, bounding the half space
${\mathbb R}^8_{+}$.
Let $A$ denote the
``arrow'' map $$A\co M_1\times M_2\longrightarrow S^7, \; A(x,y)=\frac{x-y}{|x-y|}.$$
Also consider the composition ${\phi}\co S^4_1\times S^4_2\longrightarrow S^3_1\times S^3_2\longrightarrow S^6$
of ${\Sigma} h\times {\Sigma} h$ with the arrow map in ${\mathbb R}^7$.

The $5$-cells $D^5_i$ attached to $S^4_i$ map to $M_i$, and composing with $A$ one has
a map ${\psi}\co D^5_1\times D^5_2\longrightarrow S^7$.
Pick a point $z\in S^6$ which is a regular value both of $\phi$ and $\psi$. The preimage of $z$ in $S^3_1\times S^3_2$
is a single point, and since the fiber of the Hopf map (and of its suspension) is a circle with the Lie group framing, ${\phi}^{-1}(z)$
is an Arf invariant $1$ torus $T$, the generator of
${\Omega}_2^{fr}\cong{\pi}_2^{S}$, cf \cite{Stong}. (Related stable homotopy invariants of higher dimensional links
have been considered by several authors in the context of link homotopy, see for example \cite{MR}, \cite{Kos}.)
However ${\phi}^{-1}(z)$ is the boundary
of the $3$-manifold $M^3={\psi}^{-1}(z)\subset D^5_1\times D^5_2$. Note that $T$ is codimension $6$ in
$S^4\times S^4$, while
$M^3$ is codimension $7$ in $D^5_1\times D^5_2$. Consider the pullback $\nu$ under $\psi$ of the normal line bundle
to $S^6$ in $S^7$. The $5$-cells $D^5_i$ may be assumed to be mapped near their boundary into the collar $C=S^3\times D^4\times I\subset D^8$
as a product ${\Sigma}h\times Id_{I}$. Using this product structure, one checks that there is a well-defined limit of $\nu$ at the boundary
$\partial M = T$, which amounts to a stabilization of the $6$-framing of $T$. This exhibits the Arf invariant $1$ torus $T$ as the framed boundary
of $M$, a contradiction.
\qed

Lemma \ref{Hopf lemma} easily yields a non-homological arbiter ${\mathcal A}$ on $D^8$. This step is analogous to
the corresponding part of the proof of  theorem \ref{uncountable} in dimension $4$:
consider the decomposition
$D^8=M\cup N$, where $M$ is the manifold constructed in the proof of lemma \ref{Hopf lemma},
and $N$ is its complement. Given a submanifold $P\subset D^8$ such that $M$ is isotopic to a submanifold
of $P$, set ${\mathcal A}(P)=1$. On the other hand, if $P$ is isotopic to a submanifold of $N$,
set ${\mathcal A}(P)=0$. The condition (2) in lemma \ref{Hopf lemma} ensures that this consistently defines a
partial arbiter on the ``cone'' (with respect to
the partial ordering induced by inclusion) over $M$ and the corresponding cone under $N$, see figure \ref{fig:partial ordering}.
Now pick
any arbiter on $D^8$, for example the rational homological arbiter ${\mathcal A}_{\mathbb Q}$, and use it to extend
${\mathcal A}$ to all other submanifolds of $D^8$. The condition (1) in lemma \ref{Hopf lemma} shows that
${\mathcal A}$ is different from the homological arbiters on $D^8$.

Suspending the Hopf map further, one gets the analogue of lemma \ref{Hopf lemma} and a non-homological arbiter on $D^{2n}$ for any
$n>3$. In fact, the relevant structure used in the proof of the lemma is a non-trivial square in the stable homotopy ring of spheres.
As another illustration of this construction, consider the Hopf fibration $H\co S^7\longrightarrow S^4$. Following the proof of
lemma \ref{Hopf lemma} with ${\Sigma} h$ replaced by ${\Sigma}^3 H$ (and suspending ${\Sigma}^3 H$ further), one gets arbiters on $D^{2n}, n\geq 8$.
Here at least three suspensions are needed to ensure that the dimension of the singular set $S$ is smaller than the dimension of the attaching sphere of $M$, so that \ref{Hopf lemma} (1) holds.
The following theorem summarizes the results discussed above, the case $n=3$ is proved in section \ref{D6} below.

\begin{theorem} \sl
To each non-trivial square in the stable homotopy ring of spheres
there is associated a local topological arbiter not induced by homology on
$D^{2n}$ for sufficiently large $n$. In particular, there exist non-homological local arbiters on $D^{2n}$ for each
$n\geq 3$.
\end{theorem}

\subsection{Arbiters on $\mathbf{D^6}$} \label{D6} Constructing non-homological arbiters on $D^6$ requires a bit more care.
Indeed, note that a direct analogue of the proof of lemma \ref{Hopf lemma} (with the actual Hopf map, rather than its suspension) does not hold.
(A conceptual explanation for this is provided by Lemma \ref{trivial in homology} below.) Specifically, to recall the notation consider
a standard $2$-sphere $X=S^2\subset \partial D^6$, let $C=S^2\times D^3\times I\subset D^6$ be a collar on its regular neighborhood,
and denote by $Y$ a copy of the $2$-sphere pushed into $D^6$, $Y=S^2\times\{0\}\times\{1\}$.
Let $h\co S^3\longrightarrow Y$ be the Hopf map. Consider
a general position map $f\co D^4\longrightarrow D^6$ whose restriction $f|_{\partial D^4}$ is $h$, and define $M$ to be the collar
$C$ union a regular neighborhood of the
image of $f$ in $D^6$.

Analogously to the proof of \ref{Hopf lemma}, the attaching $2$-sphere is non-trivial in homology of
$Z:=S^2\cup_{h} D^4(={\mathbb C}P^2)$, and the singular set $S$ of a general position map $f$ is still $2$-dimensional. However, in contrast with
\ref{Hopf lemma}, we show in the
lemma below that the attaching $2$-sphere is in fact trivial in $H_2(M;{\mathbb Z}/2)$ for any map $f$.
To illustrate this point, we consider a specific map
$f\co ({\mathbb C}P^2, {\mathbb C}P^1)\longrightarrow (D^4, S^3)$.  Take a Chern class one normal circle bundle $B$ over $Y$ in $D^6$.
(To be more specific, $B$ is considered as a subset of the boundary $Y\times S^3$ of a regular neighborhood of $Y$, and the fiber
of $B\longrightarrow Y$ over $y\in Y$ is the fiber $h^{-1}(y)\subset \{ y\}\times S^3$ of the Hopf fibration $h$.)

Cone $B$ to the origin of
$D^6$ and also consider the mapping cylinder of the fibration $B\longrightarrow Y$, that is the circle fiber over each
point $y\in Y$ is coned off to $y$. (This mapping cylinder  gives a map of the Chern class $1$ disk bundle into $D^6$, bounded by $B$.) This construction yields a map of $CP^2$ into $D^6$ which meets $\partial D^6$ in the $2$-sphere $X$ and the singularities may be arranged to be contained in the suspension $\Sigma$ of a single fiber $F$ of the circle bundle $B$.
The suspension is formed by the disk fiber bounded by $F$ union with ($F$ coned off
to the center of $D^6$.)
Note that this suspension $\Sigma$ represents the same class in $H_2({\mathbb C}P^2)$ as $[CP^1]$. Moreover, when $\Sigma$
is mapped into $D^6$ it is flattened into a plane, thereby killing the attaching $2$-sphere $X$ in $H_2(M;{\mathbb Z})$.

The reader may find it useful to consider a construction analogous
to the above for $({\mathbb R}P^2, {\mathbb R}P^1)$ mapped into $(D^3, S^2)$ where the Chern class $1$ bundle
$B$ is replaced by the M\"obius band. The singularity in this case is formed by a circle, generating
${\pi}_1({\mathbb R}P^2)$, flattened in a line through the origin in $D^3$.

The following lemma shows that the behavior exhibited in the example above holds (at least with ${\mathbb Z}/2$ coefficients) for {\em any}
general position map $f$.

\begin{lemma} \label{trivial in homology} \sl
Let $f\co {\mathbb C}P^2\longrightarrow D^4$ be a general position map taking ${\mathbb C}P^1$ to a $2$-sphere in $\partial D^4$.
Then the class $[{\mathbb C}P^1]$ is trivial
in $H_2(f({\mathbb C}P^2);{\mathbb Z}/2)$.
\end{lemma}

Before proving this lemma, observe that as a consequence the construction given in lemma \ref{Hopf lemma} of an arbiter
on $D^8$ fails to give a {\em non-homological} arbiter on $D^6$: since the attaching $2$-sphere is trivial in
$H_2(M;{\mathbb Z}/2)$, there is a homological obstruction to embedding two copies of $M$ into $D^6$ forming the
Hopf link in $\partial D^6$. A non-homological arbiter on $D^6$ is constructed further below in section \ref{non-homological in six dimensions}.

{\em Proof of lemma \ref{trivial in homology}.} Let $M$ be a regular neighborhood of $f({\mathbb C}P^2)$ in $D^6$. Attaching another
$6$-cell to $D^6$, consider $M$ as a submanifold of $S^6$. Suppose the statement of the lemma is false for this map $f$, then in
particular $[{\mathbb C}P^1]\neq 0\in H_2(\partial M;{\mathbb Z}/2)$. $\partial M$ is a $5$-manifold,
and there exists \cite{Thom} a $3$-manifold $N^3\subset \partial M$ representing the Poincar\'{e} dual $PD([{\mathbb C}P^1])
\in H_3(\partial M;{\mathbb Z}/2)$. Push $N$ outward to get a copy $N'$ which links the $2$-sphere
$Y$ in $S^6$. Consider the arrow map $A\co S^2\times N\longrightarrow S^5$, and composing with the
Hopf map one has $A'\co S^3\times N\longrightarrow S^5$. Since the linking number between $Y=S^2$ and $N'$ is $1$ (mod $2$), the preimage
$A'^{-1}(y)$ is a collection of an odd number of circles, each one with the Lie group framing.
We will give an argument assuming the preimage $A'^{-1}(y)$ is a single such circle $C$, the general case is proved analogously.

Since $N$ is disjoint from $M$, $A'$ extends to a map $D^4\times N\longrightarrow S^5$, and therefore
$A'^{-1}(y)$ is null-cobordant: there exists a framed (perhaps non-orientable) surface $F\subset D^4\times N^3$, $\partial F=C$.
Consider $D^4$ as the unit ball in ${\mathbb R}^4$; $C$ is a circle with non-trivial framing in $\partial D^4\times N$.
Close off $F$ to get $\overline F\subset {\mathbb R}^4\times N^3$ whose normal bundle $\nu$ is trivialized over $F$ but
$w_2({\nu})=1$. The proof of the lemma concludes with the following proposition.

\begin{proposition} \label{Whitney class} \sl
Suppose $\overline F$ is a surface embedded in $N^3\times {\mathbb R}^4$, where $N$ is a (perhaps non-orientable)
$3$-manifold. Denoting by ${\nu}$ the normal bundle over $\overline F$, the following relation holds between the
Stiefel-Whitney classes of $\nu$: if $w_1({\nu})=0$ then $w_2({\nu})=0$.
\end{proposition}

{\em Proof.} Consider the composition $g$ of the embedding $e$ with the projection onto $N$,
$g\co {\Sigma}\hookrightarrow N\times {\mathbb R}^4 \longrightarrow N$, and homotope $e$ so that $g$ is an immersion.
Denoting by ${\nu}'$ the normal bundle of the immersion $g$, note that since $\overline F$ is codimension one in $N$,
$w_1({\nu}')=0$ implies $w_2({\nu}')=0$. By naturality the same conclusion holds for $\nu$ as well. \qed

\subsection{Non-homological arbiters on $D^6$.} \label{non-homological in six dimensions}
In this section we present a modification of the argument in lemma \ref{Hopf lemma} which works
in dimension $6$. That is, we prove the existence of submanifolds $(M,S^2)\subset (D^6,S^5)$ such that

(1) \parbox[t]{14.25cm}{for any coefficient ring $R$ the attaching $2$-sphere
of $M$ is non-trivial in $H_2(M; R)$, and}

(2) \parbox[t]{14.25cm}{The components of the Hopf link $S^2\cup S^2\subset S^5$ do not bound
disjoint embeddings of two copies of $M$ into $D^6$.}

The construction of $M$ relies on the map $f\co D^4\longrightarrow D^6$ introduced at the beginning of
section \ref{D6}. Recall that the singular set $S(f)$ of that map is a $2$-sphere in $D^4$, representing
the same class as $[{\mathbb C}P^1]$ (and the same as the class of the attaching $2$-sphere $[X]$)
in the second homology of ${\mathbb C}P^2=S^2\cup_h D^4$. We will show that
a suitable $M$ may be constructed as a regular neighborhood of a map $g$: punctured $S^2\times S^2\longrightarrow D^6$.
Consider a circle linking the $2$-sphere $S(f)$ in $D^4$ and perform a $1$-surgery on $D^4$ along this circle, the result
$D^4\smallsetminus (S^1\times D^3)\cup (D^2\times S^2)$ is a punctured $S^2\times S^2$. One checks that
the $5$-dimensional $2$-handle $D^2\times D^3$ (the trace of the surgery) may be embedded in the complement of
$f(D^4)$. Denoting by $g$ the resulting map: punctured $S^2\times S^2\longrightarrow D^6$, note that the singular set of
$g$ is the same $2$-sphere, $S(g)=S(f)$.

Let $M$ be a regular neighborhood of the image of $g$, and denote by $a, b$ the generators of $H_2(S^2\times S^2)$,
corresponding respectively to the core and the co-core of the $2$-handle.
The following relation holds in $H_2(M)$: $[X]=[S(g)]+b$, where $X$ is the attaching sphere of $M$. It follows
that $[X]\neq 0\in H_2(M)$. As in the beginning of section \ref{D6}, it may be helpful to follow this
construction in the half-dimensional real context, where the analogue of $g$ maps a punctured torus into $D^3$.

The final observation is that the proof of (2) in lemma \ref{Hopf lemma} still goes through when the $5$-cells $D^5_i$ are replaced by
an arbitrary orientable punctured $5$-manifold attached to the $3$-sphere by the suspension of the Hopf map. The same argument in
dimension $6$ establishes for $M$ the property (2) above. A non-homological arbiter ${\mathcal A}$ is then defined by
setting ${\mathcal A}(N)=1$ for any $N\supset M$, ${\mathcal A}(N)=0$ for $N\subset D^6\smallsetminus M$, and extending
it to an arbiter on all submanifolds of $D^6$.

\appendix
\section{Arbiters and surgery} \label{sec: Surgery}

\begin{figure}[ht]
\includegraphics[width=4.5cm]{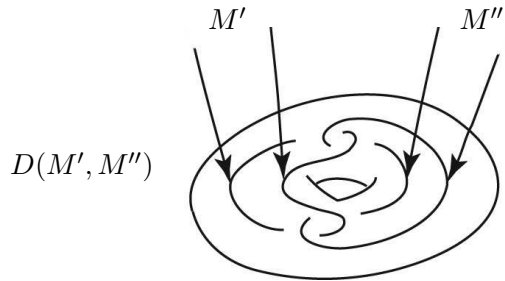}
{\small    \put(-120,95){$M'$}
    \put(-25,95){$M''$}
    \put(-195,40){$D(M',M'')$}}
    \caption{Bing Double $D(M', M'')$.}
    \label{fig:Bing double}
\end{figure}

An additional property relating the notion of a topological arbiter
to $4$-dimensional surgery (via the $(A,B)$-slice reformulation \cite{F, F1})
is stated as axiom (4) below.
As in definition \ref{def: local arbiter} of a local arbiter,
$\mathcal{M}$ denotes $\mathcal{M}_{D^4}=\{(M,{\gamma})|\, M$ is a codimension zero, smooth,
compact submanifold of $D^4$, ${\gamma}\subset
\partial M$, and $M\cap\partial D^4$ is a tubular neighborhood of an
unknotted circle ${\gamma}\subset S^3 \}$.
Given $(M',{\gamma}'), (M'',{\gamma}'') \in{\mathcal M}$,
consider their ``double'' $D(M',M'')=(S^1\times D^2\times I)\cup
(M'\cup M'')$ where $M',M''$ are attached to $S^1\times D^2\times
\{1\}$ so
that ${\gamma}',{\gamma}''$ form the Bing double of the core of the
solid torus, figure \ref{fig:Bing double}.
Let $\gamma$ denote the core of $S^1\times
D^2\times\{ 0\}$.

Note that since the components of the Bing double of the attaching curve of a $2$-handle $H=D^2\times D^2$
bound disjoint $0$-framed disks in $H$, there is an embedding $(D(M',M''),{\gamma})\subset (D^4, S^3)$ restricting to the given inclusions
$(M',{\gamma}'), (M'',{\gamma}'')\subset (D^4,S^3)$. In the axiom (4) below, $(D(M',M''),{\gamma})\in {\mathcal M}$ denotes {\em any}
embedding $(D(M',M''),{\gamma})\subset (D^4, S^3)$ restricting to the given inclusions $M', M''\subset D^4$. Consider the following condition
on a local arbiter:

\smallskip

(4) \parbox[t]{14.25cm}{Let $(M',{\gamma}'),
(M'',{\gamma}'') \in{\mathcal M}$ be such that
${\mathcal A}(M',{\gamma}')={\mathcal A}(M'',{\gamma}'')=1$. Then
${\mathcal A}(D(M',M''),{\gamma})=1$.}

\smallskip

\begin{lemma} \label{arbiter-surgery}{\sl If a local arbiter ${\mathcal A}\co {\mathcal{M}}\longrightarrow
\{0,1\}$ (definition \ref{def: local arbiter}) satisfies the axiom (4) above, it provides an obstruction to
topological $4-$dimensional surgery for the non-abelian free groups.}
\end{lemma}

{\em Proof.} We assume the reader is familiar with the $A-B$ slice formulation \cite{F, FL} of the
surgery conjecture. Let ${\mathcal A}$ be an arbiter satisfying (4).
Suppose the Borromean Rings $Bor=(l_1,l_2,l_3)$ is an $A-B$ slice
link. Consider the given decompositions $D^4=A_i\cup B_i$,
$i=1,2,3$, where the attaching curves ${\alpha}_i, {\beta}_i$ form the Hopf link in $S^3$.
Renaming the two parts of each decomposition if
necessary, using the duality axiom (3) in definition \ref{def: local arbiter}, one may assume that ${\mathcal A}(A_i, {\alpha}_i)=1$ for each
$i$. Since $Bor$ is assumed to be $A-B$ slice, there are disjoint
embeddings $f_i\co A_i\longrightarrow D^4$ such that
$f_i({\alpha}_i)=l_i$,
$i=1,2,3$. Consider the Hopf link $H=(a,b)$, and consider the first
two components $l_1, l_2$ of $Bor$ as the Bing double of the core of
a tubular neighborhood of $a$, and $l_3=b$. Viewed from this
perspective, the submanifolds $f_1(A_1), f_2(A_2)$
combine to give an embedding of a double $f\co(D(A_1,A_2),{\gamma})\longrightarrow (D^4,a)$. By axiom (4),
${\mathcal A}(D(A_1,A_2))=1$. Consider the decomposition $D^4=A\cup B$, where
$A=f(D(A_1,A_2))$, $B=D^4\smallsetminus
f(D(A_1,A_2))$. The attaching curves of $A,B$ form the Hopf link $a,b$.
Since $f_3(A_3)\subset B$, by axiom (2)
${\mathcal A}(B)=1$. Therefore ${\mathcal A}(A), {\mathcal A}(B)$ both equal $1$, a contradiction
with the axiom (3). \qed

The $A$-$B$-slice problem, a reformulation of the topological $4$-dimensional surgery conjecture, was
introduced by the first author in \cite{F, F1} and further investigated in \cite{FL}.
It is believed that the Borromean Rings are not $A$-$B$ slice, and this would imply that
surgery fails for the non-abelian free groups.
A stronger conjecture posed in these references asserted that the Borromean Rings are not even a {\em weakly
$A$-$B$ slice} link. ``Weakly A-B slice'' means that the embedding data of the submanifolds $A_i, B_i$ is not
taken into account. The second author showed in \cite{K0} that the Borromean Rings are in fact weakly $A$-$B$ slice
(section \ref{sec: embeddings} above has a more detailed discussion of this construction.) This result puts the focus
on the actual (rather than ``weak'') $A$-$B$ slice condition, where one keeps track of the embedding of the submanifolds into the
$4$-ball. The notion of a topological arbiter, additionally satisfying the condition (4) above, axiomatizes the
properties which an obstruction to surgery should satisfy in this context.
The techniques used in the construction of topological arbiters in section \ref{sec: Local arbiters}
do not give one a sufficient control over their behavior under Bing doubling. Therefore an important open question, relevant
for the surgery conjecture, is:
{\em Does there exist an arbiter ${\mathcal A}\co {\mathcal{M}}
\longrightarrow \{0,1\}$ satisfying axioms $1-4$?}

\section{Multi-arbiters.} \label{sec: multiarbiters}
This section introduces two multi-valued generalizations of topological arbiters,
and investigates their properties and applications. In particular, a non-homological multi-arbiter is constructed in
section \ref{Toda section}, and an application to percolation
theory is given in \ref{percolation section}.

Consider a $d$ dimensional manifold $M$ divided into $k$ pieces $M_i$, $M=\cup_i M_i$, with disjoint interiors, colored $1,2,\ldots, k$. To
avoid pathology for PL or smooth manifolds, we require  that each piece is a union of generic cells (dual to a triangulation).
Let $D_k$ be such a decomposition of $M$.

\begin{definition}{({\em ``Power set multiarbiters''})}  \label{power set multiarbiters def}
A {\em $k$-arbiter} on a manifold $M$ is a function from the set of
$k$-decompositions of $M$ to $2^{2^k}$,
$f\co \{ D_k\} \longrightarrow 2^{2^k}$, obeying the axioms below. For a subset $S\subset
\{ 1,\ldots, k\}$, i.e. $S\in 2^k$, let $|S|\subset M$ be $|S|=\cup_{i\in S} M_i$. Thus $|.|$ depends on $D_k$ (indicated by
a superscript below) and $f$ will be interpreted as a rule for assigning $0$ or $1$ to the submanifolds $|S|$. Let
$\overline S$ denote the complement of $S$. The following axioms are imposed on all proper subsets of $\{ 1,\ldots, k\}$, other than
the empty set.


(1) \parbox[t]{14.25cm}{$f(|S|)+f(|\overline S|)=1$.}

(2) \parbox[t]{14.25cm}{$f(|S_1|)=f(|S_2|)=1 \Rightarrow f(|S_1\cap S_2|)=1$.}

(3) \parbox[t]{14.25cm}{If $S\subseteq T$ and $f(|T|)=1$ then $f(|S|)+f(|T\smallsetminus S|)\geq 1$.}

(4) \parbox[t]{14.25cm}{If $D_k^1$ and $D_k^2$ are two $k$-decompositions so that for $i\in S$
$M_i^1\subseteq M_i^2$ then $f(|S|^1)=1 \Rightarrow f(|S|^2)=1$.}

(5) \parbox[t]{14.25cm}{If there is an isotopy of $M$ carrying $D_k^1$ to $D_k^2$ then
$f(D^1_k)=f(D^2_k)$.}
\end{definition}

The intuition is that $f(|S|)=1$ corresponds to $|S|$ carrying some predetermined homological cycle.
The definition can be applied either to closed manifolds $M$ or manifolds with boundary. In the latter case
it seems most interesting when a $k$-decomposition of $\partial M$ is fixed in advance.

To absorb the definition consider a few examples. First let $M$ be the real projective space ${\mathbb R}P^d$,
set $k=d$, and define $f(|S|)=1$ iff $|S|$ carries the non-trivial ${\mathbb Z}/2$-homology class of ${\mathbb R}P^d$ in dimension
cardinality$(S)=:c(S)$. Notice when $d=2$, one may identify this $2$-arbiter with the unique arbiter on ${\mathbb R}P^2$
described in the introduction. Similarly, picking a field $F$ yields $d$-arbiters on the $2d$- and $4d$-dimensional
manifolds ${\mathbb C}P^d$ and ${\mathbb H}P^d$.

\subsection{} \label{cube} A relative example: take a $d$-cube
$[-1,1]^d$ with a $d$-decomposition of the boundary with the $i$th piece equal to the opposite
$(d-1)$-cubes along the $i$th coordinate. Set $f(|S|)=1$ iff $$H_c(|S|, |S|\cap \partial I^d; F)
\longrightarrow H_c(I^d, |S|\cap \partial I^d; F)$$ is onto, where $c=c(S)$ and $F$ is some chosen field. Axioms 1, 2 and
3 are elementary applications of duality. Axiom 4, the ``greedy'' axiom, is also (literally) natural if one has homological
examples in mind. Finally, axiom 5 is topological invariance (the isotopy of axiom 6 is assumed to be identity on boundary
in the case of a fixed decomposition on $\partial M$.) Note that there is a unique homological arbiter (independent
of the coefficient field $F$) on the $3$-cube.

Here are a few observations on $k$-arbiters.

\subsection{} \label{Toda section}
A non-homological relative multiarbiter on the $3$-cube may be constructed as follows. When the cardinality $c(S)=1$, set $f(|S|)=0$ if
$|S|$ may be isotoped (rel boundary) in $I^3$ into the collar on $|S|\cap \partial I^3$, and $f(|S|)=1$ otherwise. Dually, for $c(S)=2$
set $f(|S|)=1$ if $|S|$ contains a $2$-disk $(D^2,\partial D^2)\subset (|S|, |S|\cap \partial I^3)$ and $f(|S|)=0$ otherwise. This definition
makes it clear that $f(|S|)+f(|\overline S|)=1$ for $c(S)=1,2$, implying axiom 1.

Axiom 2 holds since given $S_1\neq S_2\subset \{1,2,3\}$ of cardinality
$2$, if both $|S_1|, |S_2|$ contain disks then their intersection is homologically essential (rel boundary) and therefore
$|S_1\cap S_2|$ cannot be isotoped into a collar. To establish (3), suppose $T=\{a,b\}$ and $f(|T|)=1$. Then $|T|$ contains a disk $D$
as above, which is ``colored''  by two colors $a$, $b$ (this is induced by the decomposition $|T|=M_a\cup M_b$ in the notation of
definition \ref{power set multiarbiters def}). This is a relative $2$-arbiter on the disk $D$, and precisely one of the two colors, say
$a$, goes across. This implies $f(|a|)=1$, concluding the proof of (3).

A similar construction yields non-homological multi-arbiters on higher-dimensional cubes. Start with the homological $d$-arbiter on $I^d$ with
some field coefficients $F$ ($d\geq 3$), discussed in section \ref{cube},
and modify it in highest and lowest dimensions: if $c(S)=1$, set $f(|S|)=0$ if and only if
$|S|$ may be isotoped rel boundary into the collar on $|S|\cap \partial I^d$. If $c(S)=d-1$
set $f(|S|)=1$ iff $|S|$ contains a standard $(d-1)$-disk $(D^{d-1},\partial D^{d-1})\subset (|S|, |S|\cap \partial I^d)$. The axiom check is analogous to that
for the $3$-cube above.

\subsection{} \label{percolation section}
Recall the homological $k$-arbiter mentioned for the
$k$-cube in section \ref{cube}:
$f(|S|)=1$ iff
$$H_c(|S|, |S|\cap \partial I^d; {\mathbb Z}_2)
\longrightarrow H_c(I^d, |S|\cap \partial I^d; {\mathbb Z}_2)$$
is onto. It is quite easy to see (by induction on dimension) that for at least one singleton
$a\in \{ 1,2,\ldots, k\}$, $f(a)=1$. A statement in percolation follows immediately. Divide $I^d$ into
generic Voronoy cells whose centers are selected by a uniform Poisson process.
Color each cell black with probability $=1/d$. Then the probability that any two opposite
faces of $I^d$ are connected by a chain  of black cells must be $\geq 1/d$. The proof follows directly from consideration
of intersection theory and symmetry.

Consider another version of multiarbiters, whose definition is motivated by Poincar\'{e} duality and
which is a direct generalization of the $\{0,1\}$-valued arbiters (definition \ref{def: arbiter} in the introduction.)

\begin{definition}{({\em ``Poincar\'{e} multiarbiters''})} \label{def: Poincare multiarbiter}
Let $M$ be a codimension zero submanifold of a $d$-dimensional manifold $W$.
A {\em multiarbiter} $f$ associates to $M$ a number $f(M)\in\{0,\ldots,d-1\}$, subject to the axioms:

(1) \parbox[t]{14.25cm}{If $M$ is ambiently isotopic
to $M'$ in $W$ then $f(M)=f(M')$.}

(2) \parbox[t]{14.25cm}{If $M\subset M'$ then
$f(M)\leq f(M')$.}

(3) \parbox[t]{14.25cm}{Suppose $A,B$ are codimension zero submanifolds such that $W=A\cup B$, with $A\cap B
=\partial A\cap \partial B$.
Then $f(A)+f(B)=d-1$.}

(4) \parbox[t]{14.25cm}{$f(A\cap B)\geq f(A)+f(B)-d$.}
\end{definition}

In each dimension $d$, the homological multiarbiter on the projective space ${\mathbb R}P^d$ may be defined by setting
$f(M)=$maximal $k$ such that $H_k (M;{\mathbb Z}/2)$
maps onto $H_k ({\mathbb R}P^d;{\mathbb Z}/2)$. The axioms (3)-(4) follow from duality and intersection theory.

The construction in section \ref{Toda section} may be used to give an example of a non-homological Poincar\'{e} multiarbiter
on ${\mathbb R}P^3$. Given a codimension zero submanifold $M\subset {\mathbb R}P^3$, set $f(M)=2$ iff $M$ contains
a standard copy of ${\mathbb R}P^2$. Set $f(M)=0$ iff $M$ can be isotoped into a $3$-ball. $f(M)$ is defined to be $1$ on all other
submanifolds. The duality axiom (3) holds since the complement of ${\mathbb R}P^2$ is a ball, and (4) is a consequence of
intersection theory.

Recall from the introduction that if a manifold admits an open book decomposition (for example, any odd-dimensional
manifold), then it does not support a topological arbiter.
In contrast, any odd-dimensional manifold admits the ``trivial'' Poincar\'{e} multiarbiter $f$, assigning
the value $f(M)=(d-1)/2$ to each submanifold $M$. It is easy to see that this is the unique
multiarbiter (in the sense of definition \ref{def: Poincare multiarbiter}) on odd-dimensional spheres.
The classification of multiarbiters on other manifolds of dimension $\geq 3$ is an interesting open problem.


\bigskip














\begin{thebibliography}{10}

\bibitem{Dwyer} W. Dwyer, {\em Homology, Massey products and maps
between groups}, J. Pure Appl. Algebra 6 (1975), 177-190.

\bibitem{F0} M. Freedman, {\em Percolation on the projective plane}, Math. Res. Lett. 4 (1997), 889--894.

\bibitem{F} M. Freedman, {\em A geometric reformulation
of four dimensional surgery}, Topology Appl. 24 (1986), 133-141.

\bibitem{F1} M. Freedman, {\em Are the Borromean rings
$(A,B)$-slice?}, Topology Appl., 24 (1986), 143-145.

\bibitem{FL} M. Freedman and X.S. Lin, {\em On the $(A,B)$-slice
problem}, Topology Vol. 28 (1989), 91-110.

\bibitem{Hatcher} A.  Hatcher, Algebraic topology. Cambridge University Press, Cambridge, 2002.

\bibitem{Kos} U. Koschorke, {\em Link maps and the geometry of their invariants}, Manuscripta Math. 61 (1988), 383--415.

\bibitem{K0} V. Krushkal, {\em A counterexample to the strong version of Freedman's conjecture},
Ann. of Math. (2) 168 (2008), 675--693.

\bibitem{Krushkal} V. Krushkal,
{\em Link groups of $4$-manifolds},  arXiv:math/0510507.

\bibitem{Lawson} T. Lawson, {\em Open book decompositions for odd dimensional manifolds},
 Topology 17 (1978), 189-192.

\bibitem{MR}  W.S. Massey and D. Rolfsen, {\em Homotopy classification of higher-dimensional
links}, Indiana U. Math. J. 34 (1985), 375-391.

\bibitem{Milnor} J. Milnor, {\em Link Groups}, Ann. Math 59 (1954), 177-195.

\bibitem{Quinn} F. Quinn,
{\em Open book decompositions, and the bordism of automorphisms},
Topology 18 (1979), 55-73.

\bibitem{Stong} R.E. Stong, Notes On Cobordism Theory, Princeton University Press, 1968.


\bibitem{Thom} R. Thom, {\em Quelques proprietes globales de varietes differentiables}, Comment. Math. Helv.
28 (1954), 17-86.

\bibitem{Winkelnkemper}  H.E. Winkelnkemper, {\em Manifolds as open books}, Bull. Amer. Math. Soc. 79 (1973), 45-51.

\end{thebibliography}
\end{document}